\documentclass{amsart}
\usepackage{amsthm,amsmath,amsfonts,amssymb,amscd,mathrsfs,graphics}
\usepackage{latexsym}
\usepackage{graphicx}
\usepackage{rotating}

\usepackage{enumerate}

\newtheorem{thm}{Theorem}[section]

\newtheorem{lem}[thm]{Lemma}
\newtheorem{prop}[thm]{Proposition}
\newtheorem{cor}[thm]{Corollary}

\newcommand{\Aut}{\operatorname{Aut}}

\newcommand{\Id}{\operatorname{Id}}

\makeatletter
\renewcommand{\@seccntformat}[1]{\S{\csname
the#1\endcsname}\hspace{0.5em}}
\makeatother

\begin{document}

\title [Weak Cayley table groups of some crystallographic groups] {Weak Cayley table groups of  some crystallographic groups}

\author{Stephen  P.  Humphries,  Rebeca A. Paulsen}
  \address{Department of Mathematics,  Brigham Young University, Provo, 
UT 84602, U.S.A.
E-mail: steve@mathematics.byu.edu,  lakabecky@gmail.com} 
\date{}
\maketitle

\begin{abstract}  

For a group $G$, a {\it weak Cayley isomorphism} is a bijection $f:G \to G$ such that $f(g_1g_2)$ is conjugate to $ f(g_1)f(g_2)$ for all $g_1,g_2 \in G$. They form a group $\mathcal W(G)$ that is the  group of symmetries of the weak Cayley table of $G$. We determine $\mathcal W(G)$ for each of the seventeen wallpaper groups $G$, and for some other  crystallographic groups.

\medskip

\end{abstract}

\section{Introduction}

For a group $G$ the {\it weak Cayley table group} $\mathcal W(G)$ is the set of bijections $f:G \to G$ such that $f(g_1g_2)\sim f(g_1)f(g_2)$ for all $g_1,g_2 \in G$. Here $\sim$ denotes conjugacy in $G$. An element of $\mathcal W(G)$
is called a {\it weak Cayley table isomorphism} (of $G$). 

  It is easy to see that $\mathcal W(G)$ is a group that contains $\Aut(G)$ and the inverse map $\iota=\iota_G:G \to G:\iota(g)=g^{-1}$. We let $\mathcal W_0(G)=\langle \Aut(G),\iota\rangle$, the subgroup of {\it trivial } weak Cayley table isomorphisms. We note that $\mathcal W_0(G)=\Aut(G)\times \langle \iota_G\rangle.$

If $G=\{g_1=1,g_2,\dots\}$, then the group $\mathcal W(G)$ acts naturally on the {\it weak Cayley table}, this being the $|G|\times |G|$ matrix whose $ij$ entry is the conjugacy class of $g_ig_j$. It is well-known that two finite groups have the same weak Cayley table if and only if they have the same $1$- and $2$-characters, in the sense of Frobenius \cite {jms}. Here, for a character $\chi$, the corresponding $1$-character is $\chi^{(1)}=\chi$, and the $2$-character is $\chi^{(2)}:G^2 \to \mathbb C, \chi^{(2)}(x,y)=\chi(x)\chi(y)-\chi(xy).$ 

We note that the notion of $k$-character only makes sense in the situation  where $G$ is finite; however the weak Cayley table is defined for any countable group, and the group $\mathcal W(G)$ makes sense for any group.

In previous papers \cite {hu,hu2,hu3} we have determined classes of groups that have the property  $\mathcal W(G) = \mathcal W_0(G)$; these include finite symmetric groups, some free groups and free products, Coxeter groups, $PSL(2,p^n)$ and some sporadic groups. Groups with this property are termed {\it trivial}.  Elements of $\mathcal W(G)\setminus \mathcal W_0(G)$ are called {\it non-trivial}. 

In this paper we consider the seventeen wallpaper groups $G$, and determine the group $\mathcal W(G)$ explicitly; while doing so we determine which are trivial. 
We recall the standard notation and definitions for these crystallographic 
groups  \cite {ja,iv}: a {\it wallpaper group} is a discrete group of isometries of the plane  whose  subgroup, $A$,  of translations  is isomorphic to $\mathbb Z^2$. Two such are considered equivalent if they are conjugate by an affine transformation. Such a group contains only translations, rotations, reflections and glide reflections as symmetries. We use the Hermann-Mauguin notation for wallpaper groups. We note that  the automorphism group $\Aut(G)$ of each  wallpaper group $G$  has been determined in \cite {gw}, however we do not make use of their result in this paper, and so our classification will give an independent proof of their result. 

One consequence of what we do is:


\begin {thm} \label {thm1.1} Let $G$ be one of the $17$ wallpaper groups.
Then we have $\mathcal W(G)\ne \mathcal W_0(G)$ if and only if either 

\noindent (i)  $G$ is a   direct product of non-abelian groups, or 

\noindent (ii)  $G$ contains only translations and rotations, with some rotation of order greater than two. 

If $G$ has rotations of order at least four then there is a subgroup $\mathcal N(G)\triangleleft  \mathcal W(G)$, all of whose non-identity elements are non-trivial,   such that
$$\mathcal W(G)=\mathcal N(G) \rtimes \mathcal W_0(G).$$ 

\end{thm}\medskip

More detailed descriptions of the groups $\mathcal W(G)$ will be given later.
\medskip 

Let $\mathcal C_m$ denote the cyclic group of order $m$. 
    Many crystallographic groups in higher dimensions have the form $\mathbb Z^n \rtimes \mathcal C_m$; for   groups of this type  we have:
\begin {thm}\label{thmh} 
The   semi-direct product  $\mathbb Z^n \rtimes \mathcal C_2, n \ge 0,$ has trivial weak Cayley table group.
 \end{thm}

          \begin {thm}\label{thmh2} Let $p$ be an odd prime and let 
 $G=\mathbb Z^n \rtimes \mathcal C_p, n \ge 1,$ be a non-abelian  semi-direct product. Then $G$ has a non-trivial weak Cayley table map.
 \end{thm}

In \S 2 we recall some facts about elements of $\mathcal W(G)$,    describe properties of wallpaper groups, and explain the strategy for finding $\mathcal W(G)$ for each wallpaper group $G$. In the remaining sections we prove the various steps outlined in the strategy.


\section {Weak Cayley table maps and wallpaper groups}

 We collect together some results on weak Cayley table isomorphisms taken from \cite{jms}. A {\it weak Cayley table map} is a  function $\varphi:G \to H$ such that $f(g_1g_2)\sim f(g_1)f(g_2)$, for all $g_1,g_2 \in G$, where $\sim$ denotes conjugacy in $H$.
 
\begin{lem}\label{lem2.0} \cite{jms}  Let  $\varphi:G \to G $  be a weak Cayley table isomorphism and let $g \in G, N\triangleleft G$. Then,

\noindent (1)  $\varphi (e) = e $ and  $\varphi (g^{-1}) = \varphi(g)^{-1}$, where  $e $ denotes the identity element of $G$. Further, if $g^2=e$, then $\varphi(g)^2=e$. 

\noindent  (2)  $\varphi (N)\triangleleft  G$. 

\noindent  (3) 
 $\varphi(gN)=\varphi (g)\varphi (N)$.

\noindent  (4)  $\varphi$  induces a weak Cayley table map $\tilde \varphi : G/N \to  G/\varphi (N).$

  \noindent   (5) $\varphi^{-1}:G \to G$  is a weak Cayley table isomorphism.
  
  \noindent (6) If $a,b \in G$, then $a \sim b$ if and only if $\varphi(a) \sim \varphi(b)$.
 \end{lem}

We have already remarked that a wallpaper group $G$ contains a translation subgroup $A$. Also $G=A$ if and only if $G$ is abelian.
If $G=A$, then certainly $\mathcal W(G)=\mathcal W_0(G)$, so we will assume that $G$ is not abelian from now on. 
The following is well-known:

\begin{lem} \label {lem2.1} Let $G$ be a non-abelian wallpaper group. Then

(i) $[G:A]<\infty$;


(ii) $g \in G$ has a finite conjugacy class if and only if $g \in A$.
\end{lem}

Here we list  presentations  for each wallpaper group. We also,  for each of the seventeen wallpaper groups,  identify a set  $F$ of coset representatives for $G/A$:
\begin{align*}
&{\bf p1} \quad \langle x,y|(x,y)\rangle; \indent F =\{1\}. \\
&{\bf p2} \quad  \langle x,y,\rho|(x,y),x^\rho=x^{-1},y^\rho=y^{-1},\rho^2\rangle;\indent F =\langle \rho \rangle. \\
&{\bf p3}  \quad   \langle x,y,\rho | (x,y), x^\rho = x^{-1}y, y^\rho = x^{-1},   \rho^3  \rangle; \indent F =\langle \rho \rangle. \\
&{\bf p4}  \quad   \langle x,y,\rho | (x,y), x^\rho = y, y^\rho = x^{-1},   \rho^4  \rangle; \indent F =\langle \rho \rangle. \\
&{\bf p6}  \quad   \langle x,y,\rho | (x,y), x^\rho = y, y^\rho = x^{-1}y,   \rho^6  \rangle; \indent F =\langle \rho \rangle. \\
&{\bf cm}  \quad  \langle x,y,\sigma|(x,y),x^\sigma=y,y^\sigma=x, \sigma^2\rangle;\indent F =\langle \sigma \rangle. \\
&{\bf pm}  \quad  \langle x,y,\sigma|(x,y),x^\sigma=x,y^\sigma=y^{-1},\sigma^2\rangle;\indent F =\langle \sigma \rangle. \\
 &{\bf pg}  \quad  \langle x,y,\gamma|(x,y),x^\gamma=x,y^\gamma=y^{-1},\gamma^2=x\rangle; \indent F = \{1, \gamma \}. \\
  &{\bf c2mm} \quad  \langle x, y, \rho, \sigma \, \vert \, (x,y), \rho^2, \sigma^2,  x^\rho = x^{-1}, y^\rho = y^{-1}, x^\sigma = y, y^\sigma = x,  (\rho\sigma)^2 \rangle; 
  \\ & \indent F =\langle \rho, \sigma \rangle. \\
  &{\bf p2mm}  \quad  \langle x,y,\rho,\sigma|(x,y),\rho^2,\sigma^2,(\rho,\sigma),x^\rho=x^{-1},y^\rho=y^{-1},x^\sigma=x,y^\sigma=y^{-1}\rangle;\\ & \indent F =\langle \rho, \sigma \rangle. \\
   &{\bf p2mg} \quad  \langle x, y, \rho, \sigma \, \vert \, (x,y), \rho^2, \sigma^2,  x^\rho = x^{-1}, y^\rho = y^{-1}, x^\sigma = x, y^\sigma = y^{-1}, (\rho\sigma)^2 = y \rangle;  \\ &\indent F = \{1, \rho, \sigma, \rho\sigma\}. \\
&{\bf p2gg} \,\,  \langle x, y, \rho , \gamma \, \vert \, (x,y), \rho ^2, \gamma^2 = x,  x^\rho  = x^{-1}, y^\rho  = y^{-1}, x^\gamma = x, y^\gamma = y^{-1}, (\rho \gamma)^2 = y  \rangle; \\
&\indent F = \{1, \rho, \gamma, \rho \gamma \}.\\
  &{\bf p3m1}  \quad  \langle x,y,\rho,\sigma|(x,y),\rho^3,\sigma^2,(\rho\sigma)^2,x^\rho=x^{-1}y,y^\rho=x^{-1},x^\sigma=y,y^\sigma=x\rangle;\\ & \indent F =\langle \rho, \sigma \rangle. \\
  &{\bf p31m} \quad  \langle x, y, \rho, \sigma \, | \, (xy), \rho^3, \sigma^2, (\rho \sigma)^2, x^\rho = x^{-1}y, y^\rho = x^{-1}, x^\sigma = x, y^\sigma = xy^{-1} \rangle;\\ & \indent F =\langle \rho, \sigma \rangle. \\
 &{\bf p4mg}\quad  \langle x, y, \rho, \gamma \, \vert \, (x,y), \rho^4, \gamma^2 = x, x^\rho = y, y^\rho = x^{-1}, x^\gamma = x, y^\gamma = y^{-1}, (\rho\gamma)^2   \rangle; \\&\indent F = \{1, \rho, \rho^2, \rho^{-1}, \gamma, \rho \gamma, \rho^2 \gamma, \rho^{-1}\gamma\}.\\
&{\bf p4mm}\quad  \langle x, y, \rho, \sigma \, \vert \, (x,y), \rho^4, \sigma^2, x^\rho = y, y^\rho = x^{-1}, x^\sigma = x, y^\sigma = y^{-1}, (\rho \sigma)^2  \rangle;\\ & \indent F =\langle \rho, \sigma \rangle. \\
&{\bf p6m} \quad  \langle x,y,\rho, \sigma \vert (x,y), \rho^6  , \sigma^2, x^\rho = y, y^\rho = x^{-1}y,  x^\sigma = x, y^\sigma = xy^{-1}, (\rho\sigma)^2 \rangle; \\ & \indent F =\langle \rho, \sigma \rangle.
\end{align*}
In these presentations translation generators are denoted by $x,y$, rotation generators by $\rho$s, reflection generators by $\sigma$s, and glide reflection generators by $\gamma$s.

 Fundamental to what we do is the following, which is the first step in the strategy for determining $\mathcal W(G)$:

\begin{lem} If $\varphi \in \mathcal W(G)$, then $\varphi(A)=A$. In particular, for $g \in G$ we have  $\varphi(gA)=\varphi(g)A$.
\end{lem}
\noindent{\it Proof}  Since $\varphi$ preserves the size of a conjugacy class (Lemma \ref {lem2.0} (6)), the first statement  follows from Lemma \ref {lem2.1} (iii). The rest follows from Lemma \ref {lem2.0} (3). \qed\medskip

The following indicates the rest of the strategy for determining $\mathcal W(G)$ (we call these items the {\it Steps} of the proof):

\noindent (1) Show that $\varphi|_A:A \to A$ is an automorphism.

\noindent (2) Show that we can compose $\varphi$ with an element of $\mathcal W_0(G)$ so that
 the resulting element of $\mathcal W(G)$ is the identity on $A$.

\noindent (3) Show that we can then compose $\varphi$ with an element of $\mathcal W_0(G)$ so that
 the resulting element of $\mathcal W(G)$ (which we still call $\varphi$) satisfies 
$\varphi(Ag)=Ag$ for all $g\in G$. This is the same as showing that $\tilde \varphi: G/A \to G/A$ is the identity.

\noindent (4) Show that we can assume that $\varphi$ fixes each of the elements in   $F$ .

\noindent (5)  Show that for   every $t \in F$ there is some $f=f_t \in F$ such that    $\varphi(at)=a^ft$ for all $a \in A$, 

\noindent (6)  Complete the description of $\mathcal W(G)$.

\medskip

The wallpaper group $G$ acts by isometries on the Euclidean space $\mathbb E^2$; we think of a wallpaper pattern corresponding to $G$ as a subset of $\mathbb E^2$ (containing $(0,0)$). 
 In the given presentation  of $G$ we have generators $x,y$ of $A$. Further, for any $a \in A$ there is a translation of $\mathbb E^2$ by a vector $v_a \in \mathbb E^2$ (say) corresponding to $a$.
  Thus $v_x, v_y$  span a lattice $\mathfrak L\subset \mathbb E^2$. We will think of $A$ as identified with the lattice $\mathfrak L$: $a=x^iy^j\in A$ corresponds to $v_a=iv_x+jv_y\in \mathbb E$. 
  The  natural action of $G$ on $\mathbb E^2$ satisfies 
  $(v_a)g=v_{ag}$ for $a \in A, g \in G$.  
 
 Now relative to the metric on $\mathbb E^2$ we have (closed) balls $B_r$ of radius $r\ge 0$ centered at $(0,0)$. Further,  each  $f \in F$ determines an element of the orthogonal group $O(\mathbb E^2)$, where $a \mapsto a^f$, and so each conjugacy class $a^F, a \in A,$ is contained in  the boundary of   the ball $B_{|a|}$.
 
For any reflection or glide reflection $r \in G$ there is a {\it line of reflection} $L(r)\le A, L(r)=\{a \in A|a^r=a\}.$ Given $G$ we let 
$$H=H(G)=\{a \in A||a^G|\ne [G:A]\}.$$ 
It is easy to see that $H$ is the union of the lines of reflection for all reflections and  glide reflections in $G$. Since $\varphi \in \mathcal W(G)$ respects the sizes of  conjugacy classes  we have:

\begin{lem}\label {lemH} For any $\varphi \in \mathcal W(G)$ we have $\varphi(H)=H$.\qed \end{lem}

Suppose that $r \in G$ is a reflection or a glide reflection. 
Let $L^\perp (r)=\{a \in A|a^r=a^{-1}\}$.
We have:
\begin{lem} \label{lemLperp}
Suppose $r \in G$ is a reflection or a glide reflection.  Let  $\beta=r^2.$ Then for $a \in A$ we have $(ar)^2=\beta$ if and only if $a \in L^\perp (r)$.
In particular, if $\sigma \in G$ is a reflection, then $(a \sigma)^2=1 $ if and only if $a \in L^\perp (\sigma).$ 
\end{lem}
\noindent{\it Proof.} \indent We have $(ar)^2=\beta$ if and only if $ar^2 \cdot  r^{-1} ar = a\beta a^r = \beta$ if and only if $aa^r = 1$ if and only if $a \in  L^\perp (r)$.\qed\medskip


For $G$ of type {\bf p2mg} let $\beta_1=\gamma^2=1, \beta_2=(\rho\gamma)^2=y$.

For $G$ of type {\bf p2gg} let $\beta_1=\gamma^2=x, \beta_2=(\rho\gamma)^2=y$.
 
 For $G$ of type {\bf p4mg} let $\beta_1=\gamma^2=x, \beta_2=(\rho^2\gamma)^2=y$.

For $G$ of any other type let $\beta_1=\beta_2=1$.  

Let $\rho_\theta$ denote rotation by $\theta$ about $(0,0)$.

Our convention for the commutator $(g,h)$ is $g^{-1}h^{-1}gh$.
 For $f \in G$ we define  $K_f=\{(a,f):a \in A\}$.
 Since $(ab,f)=(a,f)(b,f)$ it is easy to see that $K_f$ is a subgroup of $A$.
 
\begin{lem} \label {lemK} Let $f \in F$. Then

\noindent  (i) $K_f=\langle (x,f),(y,f)\rangle.$ 

\noindent  (ii) $[G:K_f] $ is finite if and only if $f$ is a non-trivial rotation. Also $K_f=\{1\}$ if and only if $f=1$. 

\noindent  (iii)  $K_f$ is a non-trivial cyclic group if and only if $f$ is a reflection or glide reflection and if $f,f'$ are reflections or glide reflections, then $K_f=K_{f'}$ if and only if $f=f'$.

\noindent  (iv) For all $a \in A$ we have $(af)^G \cap Af = \bigcup_{f^\prime \in F^\prime} K_f(af)^{f^\prime} $ where $F^\prime = \{f^\prime \in F | (f, f^\prime) \in A\}.$  
In the cases where $G = A \rtimes F,$ this is equivalent to $(af)^G \cap Af = \bigcup_{f^\prime \in C } K_f a^{f^\prime} f$ where $C $ is the centralizer of $f$ in $F.$  

\noindent  (v)  If $f$ is a rotation,  then (iv) gives the following:
\[
(af)^G \cap Af =
\begin{cases}
  a K_ff \cup a^{\rho} K_ff  & \text{if  $G$ contains $\rho_{\pi/2}$;} \\
a K_ff \cup (a\beta_1\beta_2)^{-1} K_ff &\text{if $G$ is of type {\bf p2mm, p2mg or p2gg}}; \\
  aK_ff \cup a^\sigma K_ff & \text{if $G$ is of type {\bf c2mm}};\\
  aK_ff \cup a\beta_1\beta_1^\rho K_ff & \text{if $f=\rho_\pi$ in G of type  {\bf p4mg}};\\  
  a K_ff    &\text{otherwise.}
\end{cases}
\]
If $f$ is glide reflection or reflection, (iv) implies that
\[(af)^G \cap Af=
\begin{cases}
 a K_ff \cup (a\beta_1\beta_2)^{-1} K_ff &\text{if $G$ contains   $\rho_\pi$}; \\
  a K_ff    &\text{otherwise.}
\end{cases}
\]
  If $G$ has no glide reflections this becomes
 \[(af)^G \cap Af = 
\begin{cases}
a K_ff \cup a^{-1} K_ff &\text{if G  contains $\rho_\pi$}; \\
                a K_ff    &\text{otherwise.}
\end{cases}
\]

\end{lem} 
\noindent {\it Proof}  Since $A$ is abelian the Witt-Hall identities (\cite {mks} p. 290) show that if $a,b \in A$, then $(ab,f)=(a,f)(b,f)$, from which (i)  follows.  Now (ii) and (iii)  follow from consideration of the specific groups. 

Then  (iv) follows because conjugates of the element $af$ result from conjugating by something in $A$ or something in $F.$  Conjugating by elements in $A$ gives us $K_faf$.  If $(f^\prime, f) \notin A,$ then $(af)^{f^\prime} \notin Af.$  Thus we only need to conjugate by those $f^\prime$ that satisfy $(f^\prime, f) \in A.$ 

Most of the rest is left to the reader. We will do two cases. 
First, consider the group of type {\bf p4mm}, which contains $\rho_{\pi/2}$.  
The conjugacy class of $a\rho$ contains $(a\rho)^x = a(x,\rho^{-1})\rho$ and $(a\rho)^y = a(y,\rho^{-1})\rho$, as well as $(a \rho)^\rho = a^\rho \rho.$ 
Since $(a \rho)^\sigma \notin A\rho$ we see $(a\rho)^G \cap A\rho = 
   a K_\rho \rho \cup a^{\rho} K_\rho \rho.$ 
Likewise the conjugacy class of $a\rho^2$ contains $(a\rho^2)^x = a(x,\rho^{2})\rho^2$ and $(a\rho^2)^y = a(y,\rho^{2})\rho^2$, as well as $(a \rho^2)^\rho = a^\rho \rho^2.$ 
Since $(a \rho^2)^\sigma \in aK_{\rho^2} \rho^2$ we have $(a\rho^2)^G \cap A\rho^2 = 
   a K_{\rho^2} \rho^2 \cup a^{\rho^2} K_{\rho^2} \rho^2,$ proving the first line of (v).

For a second example, suppose $\gamma$ is a glide reflection in a group with a generator $\rho =\rho_\pi.$  
Recall that we defined $\beta_1 = \gamma^2, \beta_2=(\rho \gamma)^2$. 
Then  $(a\gamma)^G$ will include 
$(a\gamma)^x = a(x,\gamma^{-1})\gamma $ and $(a\gamma)^y = a(y,\gamma^{-1})\gamma,$ thus it includes $aK_\gamma \gamma.$  
Note that $(a,\gamma) \in K_\gamma$ implies that $(a\gamma)^\gamma \in aK_\gamma\gamma.$
Lastly we consider $(a\gamma)^\rho.$ 
Since $\gamma^\rho = \beta_1^{-1}\beta_2\gamma,$ we have $
(a\gamma)^\rho = (a \beta_1 \beta_2^{-1})^{-1}\gamma.$
 But since $\beta_2^2 \in K_\gamma$ we see  that $(a\beta_1\beta_2)^{-1}K_\gamma \gamma \subset (a\gamma)^G. $
Thus $(a \gamma)^G = aK_\gamma \gamma \cup (a\beta_1\beta_2)^{-1}K_\gamma \gamma$. 
\qed\medskip

\noindent {\bf Some Automorphisms   of wallpaper groups}

Here we list some automorphism that we will use in what follows. 
We do not list inner automorphisms, which  will be denoted $I_g, g \in G$.

\noindent {\bf p3, p4, p6} 
$\psi_x: x \mapsto x, y \mapsto y,  r  \mapsto xr;
\psi_y:x \mapsto x, y \mapsto y,   r  \mapsto yr.$

\noindent {\bf cm} $\psi: (x,y,\sigma) \mapsto (x^{-1},y^{-1},\sigma)$.

\noindent {\bf pm} $\psi: (x,y,\sigma) \mapsto (x^{-1},y^{-1},\sigma)$.

\noindent {\bf pg} $\psi_y: (x,y,\gamma) \mapsto (x,y,y\gamma)$.  Also $\psi: (x,y,\gamma) \mapsto (x^{-1},y^{-1},\gamma^{-1})$.

\noindent {\bf c2mm} $\psi_{u,v,i,j}:(x,y,\rho,\sigma) \mapsto (x,y,x^uy^v\rho,x^iy^j\sigma)$ where $u+j=i+v$; also $\psi:(x,y,\rho,\sigma) \mapsto (x^{-1}, y^{-1}, \rho, \sigma).$

\noindent {\bf p2mm} $\psi_{u,v}:(x,y,\rho,\sigma) \mapsto (x,y,x^uy^v\rho,x^uy^v\sigma).$
Also $\psi:(x,y,\rho,\sigma) \mapsto (y,x,\rho,\rho\sigma)$.

\noindent {\bf p2mg} $\psi_x:(x,y,\rho,\sigma) \mapsto (x,y,x\rho,\sigma)$; also  $\psi_y:(x,y,\rho,\sigma) \mapsto (x,y,y\rho,y\sigma)$; and
$\psi:(x,y,\rho,\sigma) \mapsto (x^{-1},y^{-1},y^{-1}\rho,\sigma).$

\noindent {\bf p2gg} $\psi_{x}:(x,y,\rho,\gamma) \mapsto (x,y,x\rho,\gamma)$ and $ \psi_{y}:(x,y,\rho,\gamma) \mapsto (x,y,y\rho,y\gamma). $ Also $\psi:(x,y,\rho,\gamma) \mapsto (x^{-1}, y^{-1}, \rho, x^{-1}y\gamma).$  

\noindent {\bf p3m1} $\psi_1:(x,y,\rho,\sigma) \mapsto (x,y,y^{-1}\rho,\sigma).$

\noindent {\bf p4mg} $\psi_1:(x,y,\rho,\gamma) \mapsto (x,y,x\rho,y^{-1}\gamma);$ $\psi_2:(x,y,\rho,\gamma) \mapsto (x^{-1},y^{-1},y^{-1}\rho,x^{-1}\gamma).$ 

\noindent {\bf p4mm} $\varphi_1:(x,y,\rho,\sigma) \mapsto (x,y,y\rho,y\sigma)$.

\medskip 

\noindent {\bf Some  non-trivial weak Cayley table  isomorphisms  of wallpaper groups}

\noindent {\bf p2mm} Define $\tau$ by   $ \tau:\Bigl \{ \begin{matrix}  g \mapsto g \,\,\text{ for }  g \in A \cup A\rho\sigma; \\ g \mapsto g^\sigma  \,\, \text{ for }  g \in A\rho \cup A\sigma\end{matrix}$.

\noindent {\bf p3} Define $ \tau:\Bigl \{ \begin{matrix}  g \mapsto g \,\,\text{ for }  g \notin A\rho \cup A\rho^2; \\ g \mapsto g^\rho\,\, \text{ for }  g \in A\rho \cup A\rho^2\end{matrix}$.
\medskip

\noindent {\bf p4} For $h \in \{x,y,\rho^2\}$ define $ \tau_h:\Bigl \{ \begin{matrix} g \mapsto g &\text{ for }g \notin A\rho^2;\\ g \mapsto g^h& \text{ for } g \in A\rho^2\end{matrix}$.

For $h \in \{x,y,\rho\}$ define $ \mu_h:\Bigl  \{ \begin{matrix}  g \mapsto g &\text { for } g \in A \cup A\rho^2; \\g \mapsto g^h&  \text{ for } g \notin A \cup A\rho^2\end{matrix}$.
\medskip 

\noindent {\bf p6} For $h \in \{xy,xy^{-2},\rho^2\}  $ define $ \tau_h:\Bigl  \{ \begin{matrix} g \mapsto g &\text{ for } g \in A\rho \cup A\rho^3\cup A\rho^5;\\ g \mapsto g^h& \text{ for } g \notin A\rho \cup A\rho^3\cup A\rho^5\end{matrix}$.

For $h \in \{x^2,y^2,\rho^3\}$ define $ \mu_h:\Bigl  \{ \begin{matrix}g \mapsto g & \text{ for } g \notin A \cup A\rho^3; \\ g \mapsto g^h& \text{ for } g \in A \cup A\rho^3 \end{matrix}$.

\section {Steps (1), (2)}

Let $\varphi \in \mathcal W(G)$ where $G$ is a wallpaper group. 
In this section we will denote $A$ additively, so that the action of $g \in G$  on $A$ by conjugation is denoted $(a)g, a \in A$.
Then for all $a,b \in A$ there is some  $g=g_{a,b} \in F$ such that
\begin {align}  \varphi(a+b)=(\varphi(a)+\varphi(b))g_{a,b}.\tag   {3.1}\label {e31} \end {align} 
Thus
\begin{align}\notag 
\varphi(b)&=\varphi(a+b-a)=(\varphi(a+b)+\varphi(-a))g_{a+b,-a}\\&=(\varphi(a+b)-\varphi(a))g_{a+b,-a}.\tag {3.2} \label {e32}
\end{align}

\begin {lem} \label {lem31} For all $a,b \in A$ there is $f \in F$ such that
\begin{align} \varphi(a+b)=\varphi(a)+ \varphi(b)f.\notag 
\end{align}
\end{lem}
\noindent {\it Proof} From equation (\ref {e32}) we see that $f=g_{a+b,-a}^{-1}$ will work.\qed\medskip

Let $a,b \in A$. We wish to show that $\varphi(a+b)=\varphi(a)+ \varphi(b)$.
From equation (\ref{e31}) and  Lemma \ref {lem31}  we see that
\begin{align}\tag {3.3} \label {e33}  \varphi(a+b)\in (\varphi(a)+\varphi(b))F \cap (\varphi(a)+\varphi(b)F) \cap (\varphi(a)F+\varphi(b)).\end{align}

Since $F$ acts by orthogonal matrices we see that 

\noindent (i) any point of the form
$(\varphi(a)+\varphi(b))F$ is on the  circle $C_1$ of radius $|\varphi(a)+\varphi(b)|$ centered at $(0,0)$; 

\noindent (ii)  any point of the form
$\varphi(a)+\varphi(b)F$ is on the  circle  $C_2$  of radius $|\varphi(b)|$ centered at $\varphi(a)$;

\noindent (iii)  any point of the form
$\varphi(a)F+\varphi(b)$ is on the circle  $C_3$ of radius $|\varphi(a)|$ centered at $\varphi(b)$.

The intersection $C_1 \cap C_2 \cap C_3$ certainly contains the point $p_1=\varphi(a)+\varphi(b)$, and any other possibility for $\varphi(a+b)$.
 If $p_1$ is the only such point of $C_1 \cap C_2 \cap C_3$, then by equation (\ref {e33}) we must have $\varphi(a+b)=\varphi(a)+\varphi(b)$, and we are done.

 So now assume that $C_1 \cap C_2 \cap C_3$ also contains $p_2 \ne p_1$; then the arc $p_1p_2$ has a perpendicular bisector, $L$ say,
that contains the centers of each of $C_1, C_2$ and $C_3$. Thus
$\varphi(a), \varphi(b)$ and $(0,0)$ are on this line, and so $p_1=\varphi(a)+\varphi(b)$ is also on  $L$.


But $p_1$ was on a line  perpendicular to $L$ that also contains $p_2$, where the distance from $p_1$ to $L$ was the same as the distance of $p_2$ to $L$. It follows that $p_1=p_2$, and we have a contradiction.
Thus $\varphi|_A$ is a homomorphism; however $\varphi$ is a bijection, and so we have shown Step (1): 

\begin{lem} \label {lemhom} For any $\varphi \in \mathcal W(G)$ we have $\varphi|_A:A \to A$ is an automorphism.\qed \end{lem}

 So we may now assume that $\varphi|_A:A \to A$ is an automorphism.
   Relative to the basis $v_x, v_y$ the matrix representing  $\varphi|_A$ is an integer matrix, $M(\varphi)$ say, and since it is invertible $M(\varphi)$  must have determinant $\pm 1$. We first show:

 \begin {lem} \label {finorder} If $G$ does not have type {\bf p1,p2}, then   $\varphi|_A$ has finite order. 
 \end{lem}
 \noindent {\it Proof} 
 First assume that there is some $a \in A$ such that $(a)F$ does not span a cyclic subgroup of $A$. Let $a_1, a_2 \in (a)F$ generate a $\mathbb Z^2$ subgroup.
 If $\varphi|_A$ has infinite order, then for every $N>0$ there is $k\in \mathbb N$ such that some point of $\varphi^k((a)F)$ is outside $B_N$. 
 But all the elements of $(a)F$ are conjugate, and so all the elements of
 $\varphi^k((a)F)$ are conjugate by Lemma \ref {lem2.0}.
 Thus if one of the points of $\varphi^k((a)F)$ 
  is outside $B_N$, then they are all outside of $B_N$. 
  Thus $\varphi(a_1), \varphi(a_2)$ are both on the boundary of the ball $B_{|\varphi(a_1)|}\supset  B_N$. Now the triangle with vertices $(0,0), a_1, a_2$ is sent to the triangle with vertices $(0,0), \varphi(a_1),\varphi(a_2)$, however (as  the distance between $\varphi(a_1),\varphi(a_2)$ is at least $1$, and we can choose $N$ arbitrarily large)  the latter triangle clearly has  larger area than the former, contradicting the fact that $M(\varphi)$ has determinant $\pm 1$. 
 
 Now one can  check (i) and (ii) of the following result:
 
 \begin {lem} \label{lemnotp2}  (i) Every wallpaper group $G$ other then {\bf p1} and {\bf p2} has some $a \in A$ where $(a)F$ does not span a cyclic group. 

(ii) If $\tau$ is some isometry of the wallpaper paper pattern in $\mathbb E^2$ corresponding to the wallpaper group $G$ and $\tau((0,0))=(0,0)$, then there is some element of $ \mathcal W_0(G)$ whose action on $A$  is equal to that of $\tau$ on $\mathbb E^2$. \qed 
\end{lem} 
 
 Lemma \ref {finorder} now follows from Lemma \ref {lemnotp2} (i).\qed\medskip

 Thus by Lemma \ref {finorder}  we see that if $G$ is not of type {\bf p1,p2}, then $M(\varphi)$ is an integer matrix of finite order.
  This order is well-known to be $1,2,3,4$ or $6$, and the element is an orthogonal rotation or  a reflection. This matrix also preserves the lattice $\mathfrak L$.    It thus 
  corresponds to the action of some element of $\mathcal W_0(G)$, by Lemma \ref {lemnotp2} (ii). Composing $\varphi$ with the inverse of this element gives a new $\varphi$ where $\varphi|_A=\Id_A$, as required for Step (2).
 
 Thus we now assume that $G$ has type {\bf p2}. In this case it is easy to see that whenever $a,b,c,d \in \mathbb Z$ satisfy $ad-bc=\pm1$, then the assignments
 $$x \mapsto x^ay^b,\quad y \mapsto x^cy^d, \quad r \mapsto r,$$
 determine  an automorphism of $G$. Thus we can compose $\varphi$ with some element of $\mathcal W_0(G)$ to obtain a new $\varphi$ satisfying $\varphi|_A=\Id_A$. This 
 gives Step (2):
 \begin{prop} For $\varphi \in \mathcal W(G)$  there is   $\varphi' \in \mathcal W_0(G)$ with  $(\varphi'\circ \varphi)|_A=\Id_A$.\qed \end{prop}
 
  \section {Step (3)} 
 
 Since $A \triangleleft G$  we see from Lemma \ref{lem2.0} that $\varphi$ induces a weak
Cayley table isomorphism $\tilde \varphi:G/A \to G/A$.  Now $G/A$ is one of the groups $\mathcal C_1,\mathcal C_2, \mathcal C_3, \mathcal C_4,\mathcal C_6,
D_4= \mathcal C_2 \times \mathcal C_2, D_6,D_8,D_{12}$. It is known (\cite {hu}) that each of these groups has trivial weak Cayley table group. 
 
 \begin{lem} \label{leminvA} If $G$ is a wallpaper group having no glide reflections in $F$, with presentation as given in \S 2 with  generators $x,y,\rho,\sigma$, then the assignments 
 $$x \mapsto x^{-1}, y \mapsto y^{-1},  \rho \mapsto \rho,\sigma \mapsto \sigma,$$ determine an automorphism  of $G$ that we denote by $\iota_A$.
\end{lem} 
 \noindent {\it Proof} This follows from the fact that each relation involving $x,y$ has the form $(x,y)$ or $x^\rho=w=w(x,y)\in A$, this latter type being equivalent to
  $(x^{-1})^\rho=w(x,y)^{-1}=w(x^{-1},y^{-1})$.\qed\medskip

 Since  $\tilde \varphi:G/A \to G/A$ is an automorphism or antiautomorphism, we see that if $\rho \in G$ is a rotation that generates the subgroup of rotations fixing $(0,0)$, then we must have 
 $\tilde \varphi(A\rho)=A\rho^{\pm 1}$, and so $\varphi(\rho)\in A\rho^{\pm 1}$. If $\varphi(\rho)\in A\rho^{-1}$, then we compose $\varphi$ with $\iota$ and then with $\iota_A$, to give a new $\varphi$ which is still the identity on $A$, but which satisfies $\varphi(\rho)\in A\rho$.
 
 If $|G/A|=2$, then   $\tilde \varphi$ is the identity. So assume that $|G/A|>2$.
 
If $G$ does not have reflections or glide reflections,  then the above shows that  $\tilde \varphi$ is the identity. 

So now we assume that $G$ has reflections, but does not have  glide reflections. Then we have $G/A \cong D_4=\mathcal C_2\times \mathcal C_2,  D_6, D_8, D_{12}$. In particular, $G$ has at least two reflection cosets:
there are $\sigma_1,\sigma_2 \in G, A\sigma_1 \ne A\sigma_2,$ where $\sigma_i,i=1,2,$  are reflections. Assume  that $\varphi(\sigma_1)=a\sigma_2$.  Since $\sigma_1$ is a reflection, then $a \in  L^\perp (\sigma_2)$ by Lemma \ref {lemLperp} and Lemma \ref {lem2.0} (1). Also, $L^\perp (\sigma_1) \ne L^\perp (\sigma_2)$, since $A\sigma_1 \ne A\sigma_2$.
 Then there is $b \in  L^\perp (\sigma_1)$ such that $ab \notin  L^\perp (\sigma_2)$. But $b \sigma_1$ has order $2$, however    $\varphi(b\sigma_1)\sim \varphi(b)\varphi(\sigma_1)=ba\sigma_2$; since $ab \notin  L^\perp (\sigma_2)$ we see that $ab\sigma_2$ does not have order $2$, a contradiction.
 Thus we have proved Step (3) for groups without glide reflections:
 \begin{lem} If $G$ does not have glide reflections in $F$, then we can assume (by composing with a trivial weak Cayley table map if necessary) 
 that $\varphi|_A=\Id_A$ and that
 $\tilde \varphi:G/A\to G/A$ is the identity.\qed \end{lem} 
 
 So now assume that $G$ has a glide reflection in $F.$ There are four wallpaper groups of this type which we now consider individually.
 
 If $G$ is of type {\bf pg}, then $|G/A|=2$ and so we certainly have $\tilde \varphi=\Id_{G/A}$. 
 
 If $G$ is of type {\bf p2mg}, then $G/A\cong \mathcal C_2^2$, and the  cosets are $A, A\rho, A\sigma, A\rho\sigma$. Here $A\rho$ consists entirely of rotations;
if  $x^iy^j\sigma \in A\sigma$, then $(x^iy^j\sigma)^2=x^{2i}$; and  
 if  $x^iy^j\rho \sigma \in A\rho\sigma$, then $(x^iy^j\rho\sigma)^2=y^{2j+1}$. This shows that $\varphi(A\rho)=A\rho$. If we have $\varphi(\sigma)=x^iy^j\rho\sigma$, then $\sigma^2=1$ gives $(x^iy^j\rho\sigma)^2=1,$ so that $y^{2j+1}=1$, a contradiction. Thus we have $\tilde \varphi=\Id_{G/A}$.

 If $G$ is of type {\bf p2gg}, then $G/A\cong \mathcal C_2^2,$ and the  cosets are $A, A\rho, A\gamma, A\rho\gamma$.  Again $A\rho$ consists entirely of order $2$ rotations, while each element of $A\gamma, A\rho\gamma$  is a glide reflection and so not of order $2$. Thus we  have $\varphi(A\rho)=A\rho$. If  $\varphi(\gamma)=
 x^iy^j\rho\gamma$, then
\begin{align*}
x&=\varphi(x)=\varphi(\gamma^2)\sim (x^iy^j\rho\gamma)^2=x^iy^j\rho\gamma x^iy^j\rho\gamma\\&
=x^iy^j\rho\gamma\rho x^{-i}y^{-j}\gamma
=x^iy^j y \gamma^{-1} x^{-i}y^{-j}\gamma=y^{2j+1}.
 \end{align*} Thus $ y^{2j+1}\in x^G= \{x,x^{-1}\}$,  a contradiction. Thus we have $\tilde \varphi=\Id_{G/A}$.
 
 If $G$ is of type {\bf p4mg}, then $G/A\cong   D_8$. We describe  $\mathcal W(G/A)$.
 Order the cosets of $G/A $ as
 $$A,A\rho^2,A\rho^2\gamma, A\gamma, A\rho\gamma, A\rho^3\gamma, A\rho^3, A\rho.$$
 Then relative to this ordering the action of $\mathcal W(G/A)$ is given by the permutation group 
 $\langle (3,5,4,6),(3,4),(7,8)\rangle\cong D_8 \times \mathcal C_2.$

  Thus  the possibilities for $\tilde \varphi (A\rho)$ are $A\rho, A\rho^3$; we also see that
 $\varphi(A\rho^2)=A\rho^2$, and that  $\{A\gamma,A\rho\gamma,A\rho^2\gamma,A\rho^3\gamma\}$ are permuted by $\tilde \varphi$. If $\varphi(A\rho)=A\rho^3$, then composing with $\iota$ and then $\psi_2$ shows that we can assume $\varphi |_A=\Id_A, \tilde \varphi(A\rho)=A\rho, \tilde \varphi(A\rho^2)=A\rho^2, \tilde \varphi(A\rho^3)=A\rho^3$. 
 
Considering the other cosets we have:
\begin{align*}
&(x^iy^j\gamma)^2=x^{2i+1};\qquad (x^iy^j\rho\gamma)^2=x^{i-j}y^{j-i};\\
& (x^iy^j\rho^2\gamma)^2=y^{2j+1};\qquad  (x^iy^j\rho^3\gamma)^2= (xy)^{i+j+1}.
\end{align*}

This shows that $A\rho\gamma, A\rho^3\gamma$ are the only cosets that contain involutions, so that we must have
$\{\varphi(A\rho\gamma),\varphi(A\rho^3\gamma)\}=\{A\rho\gamma,A\rho^3\gamma\}$. 
This shows that the action of $\tilde \varphi$ is given by some element of the group $\langle (3,4),(5,6)\rangle.$
If $\tilde \varphi \ne 1$, then we must  either have 
 $\varphi(A\gamma)=A\rho^2\gamma$ or  $\varphi(A\rho \gamma)=A\rho^3\gamma$.

\noindent {\bf CASE 1}:  We first show that we cannot have 
$\varphi(A\gamma)=A\rho^2\gamma$. So assume that this is the case, and  that $\varphi(\gamma)=x^uy^v\rho^2\gamma$. Since $\gamma^2=x$ we have
\begin{align*} x&=\varphi(x)=\varphi(\gamma^2) \sim
(x^uy^v\rho^2\gamma)^2
=x^uy^v\rho^2\gamma\cdot  x^uy^v\rho^2\gamma\\
&=x^uy^v(\rho^2\gamma \rho^2) x^{-u}y^{-v}\gamma=
x^uy^v (y \gamma^{-1})  x^{-u}y^{-v}\gamma=x^uy^v y   x^{-u}y^{v}
=y^{2v+1}.
\end{align*}
Since $x^G=\{x,y,x^{-1},y^{-1}\}$ we see that $2v+1 \in \{1,-1\}$, so that $v\in \{0,-1\}$.


The next  result follows from the presentation for $G$ of type  {\bf p4mg}:
\begin{lem} \label{lembr2s} Let $G$ have type {\bf p4mg}. Let $X=\langle x^2\rangle, Y=\langle y^2\rangle$. For $a \in A$ we have
$$(a\rho^2\gamma)^G=yYa^{\rho^3}\gamma \cup x^{-1}Ya^\rho\gamma \cup Xa\rho^2\gamma\cup xy^{-1}Xa^\gamma \rho^2\gamma.\qed
$$
\end{lem}

Now $(x\gamma)^2=x^3$, giving  $x^3=\varphi(x^3)=\varphi((x\gamma)^2) \sim \varphi(x\gamma)^2\sim (\varphi(x)\varphi(\gamma))^2$, so that
$x^3 \sim (x^{u+1}y^v\rho^2\gamma)^2$.   Lemma \ref {lembr2s}  shows that 
$(x^{u+1}y^v\rho^2\gamma)^G$ has one of the following forms, for some $k \in \mathbb Z$: 
\begin{align*}
&(i) \, y^{2k+1}y^{-u-1}x^v\gamma;\qquad  
(ii)\,  x^{-1}y^{2k}y^{u+1}x^{-v}\gamma;\\
& (iii)\, x^{2k}x^uy^v\rho^2\gamma;\qquad 
(iv) \, x^{2k+1}y^{-1}x^{u+1}y^{-v}\rho^2\gamma.\end{align*}
If we have (i), then  
\begin{align*}x^3&\sim 
(y^{2k+1}y^{-u-1}x^v\gamma)^2=
y^{2k-u}x^v\gamma\cdot y^{2k-u}x^v\gamma=y^{2k-u}x^v\gamma^2y^{-2k+u}x^v
=x^{2v+1},
\end{align*}
which is a contradiction since $v \in \{0,-1\}$. 

If we have (ii), then  $x^3 \sim (y^{2k+u+1}x^{-v-1}\gamma)^2$, where
\begin{align*}
(y^{2k+u+1}x^{-v-1}\gamma)^2&=y^{2k+u+1}x^{-v-1}\gamma \cdot y^{2k+u+1}x^{-v-1}\gamma\\&=
y^{2k+u+1}x^{-v-1}\gamma^2 \cdot y^{-2k-u-1}x^{-v-1}=x^{-2v-1},
\end{align*}
which is a contradiction since $v \in \{0,-1\}$. 

If we have (iii), then  $x^3 \sim (x^{2k+u}y^{v}\rho^2\gamma)^2$, where
\begin{align*}
(x^{2k+u}y^{v}\rho^2\gamma)^2&=x^{2k+u}y^{v}\rho^2\gamma \cdot x^{2k+u}y^{v}\rho^2\gamma
=x^{2k+u}y^{v}\rho^2\gamma\rho^2   x^{-2k-u}y^{-v}\gamma\\
&=x^{2k+u}y^{v} y \gamma^{-1}  x^{-2k-u}y^{-v}\gamma=y^{2v+1},
\end{align*}
which is a contradiction since $v \in \{0,-1\}$. 

If we have (iv), then  $x^3 \sim (x^{2k+2+u}y^{v}\rho^2\gamma)^2$, where
\begin{align*} (x^{2k+2+u}y^{v}\rho^2\gamma)^2&=x^{2k+2+u}y^{v}\rho^2\gamma\cdot x^{2k+2+u}y^{v}\rho^2\gamma=
x^{2k+2+u}y^{v}\rho^2\gamma \rho^2      x^{-2k-2-u}y^{-v}\gamma\\
&=x^{2k+2+u}y^{v} y\gamma^{-1}   x^{-2k-2-u}y^{-v}\gamma=y^{2v+1},
\end{align*}
which is a contradiction since $v \in \{0,-1\}$.  
Thus CASE 1 does not happen. \medskip

\noindent {\bf CASE 2}: Here we assume that $\varphi(A\rho \gamma)=A\rho^3\gamma$ and write  $\varphi(\rho \gamma)=b\rho^3\gamma, b=x^uy^v$. 
Since $(\rho\gamma)^2=1$ we must have $(b\rho^3\gamma)^2=1$. But
\begin{align*}
(b\rho^3\gamma)^2&=x^uy^v\rho^3\gamma \cdot x^uy^v\rho^3\gamma=x^uy^v\rho^3\gamma \rho^3y^{-u}x^v\sigma=
x^uy^vy\gamma y^{-u}x^v\sigma\\
&=x^uy^{v+1}s^2y^ux^v=(xy)^{u+v+1},
\end{align*} so that $u+v+1=0$.

Now $(xy\rho\gamma)^2=1$, and so
$1=\varphi((xy\rho\gamma)^2)\sim \varphi(xy\rho\gamma)^2$, so that $\varphi(xy\rho\gamma)^2=1$. But $\varphi(xy\rho\gamma)\sim xy\cdot x^uy^v\rho^3\sigma$, and so we must have $(x^{u+1}y^{v+1}\rho^3\sigma)^2=1$. One finds 
 that
$(x^{u+1}y^{v+1}\rho^3\sigma)^2=(xy)^{u+v+3},$ so that $u+v+3=0,$ a contradiction.

Thus CASE 2 does not happen  either,  and we have proved Step (3).

\section {Step (4)} 

We need to show that we can assume (after possibly composing $\varphi$ with an element of $\mathcal W_0(G)$ or a known non-trivial weak Cayley table isomorphism)  that $\varphi$ fixes each   element of 
$F$. We do this in such a way as to preserve the fact that $\varphi|_A$ is the identity on $A$, and that $\tilde \varphi$ is the identity map on cosets of $A$. We consider 
 each type of wallpaper group.\medskip 

\noindent   {\bf p2, p3, p4, p6} cases:  Here we  have $\varphi(\rho)=a\rho, a \in A$; however the assignment $(x,y,\rho) \mapsto (x,y,a\rho)$ determines  an automorphism of $G$ in this case. To see this we note that, by the classification of isometry types, $a\rho$ is either a rotation, a translation, a reflection or a glide reflection, the last three options not being possible. 
We thus compose $\varphi$ with the inverse of this automorphism to obtain $\varphi(\rho)=\rho, \varphi(\rho^{-1})=\rho^{-1}$.
This then does the cases {\bf p2, p3}. 

For {\bf p4} we note that $\varphi(\rho^2)\sim \varphi(\rho)^2=\rho^2$, which shows that $\varphi(\rho^2)=(\rho^2)^a, a=x^iy^j \in A$. But then composing $\varphi$ with  
$\tau_x^{-i}\tau_y^{-j}$  will then give $\varphi(\rho^i)=\rho^i,0\le i<4$. 

For {\bf p6} we note that $\varphi(\rho^2)\sim \varphi(\rho)^2=\rho^2$, which shows that $\varphi(\rho^2)=(\rho^2)^a, a=x^iy^j \in A$. Thus $\varphi(\rho^2)=(xy)^k(xy^{-2})^m\rho^2$ for some $k,m \in \mathbb Z$.  But now we can compose $\varphi$ with $\tau_{xy}^{-k}\tau_{xy^{-2}}^{
-m}$ so as to be able to assume $\varphi(\rho^2)=\rho^2,
\varphi(\rho^{-2})=\rho^{-2}$.  Now  $\varphi(\rho^3)\sim \varphi(\rho)\varphi(\rho^2)=\rho^3$, so that $\varphi(\rho)=x^{2i}y^{2j}\rho^3, i,j \in \mathbb Z$. Now one can use 
$\mu_{x^2}^{-i}\mu_{y^2}^{-j}$ so as to be able to assume $\varphi(\rho^3)=\rho^3$. This does the {\bf p6} case. 
\medskip

\noindent   {\bf cm, pm} cases: Here we must have $\varphi(\sigma)=a\sigma, a \in A$, where $(a\sigma)^2=1$ since $\sigma^2=1$; however the assignment $(x,y,\sigma) \mapsto (x,y,a\sigma)$ determines  an automorphism of $G$ in these cases, and we are done.
\medskip 

\noindent   {\bf pg} case: Here we note that $x=\gamma^2$ is central in $G$, so that we must have $\varphi(\gamma)^2=x$; one finds that if $\varphi(\gamma)=x^uy^v\gamma$, then $u=0$. Thus we can compose with  a power of the automorphism $\psi_y$ of $G$ so as to get $\varphi|_A=Id$ and $\varphi(\gamma)=\gamma$. 
\medskip

We will find the following useful in what follows:
\begin{lem}\label{lemcfg} Suppose that $\varphi|_A=\Id_A$ and that $\varphi(\rho)=\rho$ for a rotation $\rho$. 
Suppose that $\sigma,\sigma\rho$ are distinct reflections in $G$ with $\varphi(A\sigma)=A\sigma$. Then $\varphi(\sigma)=\sigma$.
\end{lem}
\noindent {\it Proof} Since  $\varphi(A\sigma)=A\sigma$ we can write $\varphi(\sigma)=a\sigma, a \in A$. Since $\sigma^2=1$ we must have $a \in L^\perp (\sigma)$.
Now 
$\varphi(\rho\cdot \sigma)\sim \rho a \sigma=a^{\rho^{-1}}\rho\sigma$ and so $a^{\rho^{-1}} \in L^\perp (\rho \sigma)$. But then  $a \in L^\perp(\sigma\rho)$. Since $\sigma \ne \sigma\rho$ we have $a \in  L^\perp (\sigma) \cap L^\perp(\sigma\rho)=\{1\}$, and we are done.\qed\medskip 

%

\begin{lem}\label{lem_phi(gamma)}
Let $\varphi|_A = Id|_A$ and assume $\varphi(A\gamma) = A\gamma$ for some glide reflection coset $A\gamma.$    Then $\varphi(\gamma) = a\gamma$ implies that $a \in L^\perp(\gamma).$  
\end{lem}
\noindent {\it Proof.} \indent
Let $\beta = \gamma^2,$ so that $\beta \in L(\gamma), \beta\ne 1$.
Squaring both sides of   $\varphi(\gamma) = a\gamma$ gives us $\beta \sim  \beta aa^\gamma.$ Since $\beta$ and $aa^\gamma$ both lie on $L(\gamma),$ this implies that $\beta aa^\gamma = \beta^{\pm1}.$  
Then we have two possibilities:  either $aa^\gamma = 1$ which implies that $a \in L^\perp(\gamma),$ or $aa^\gamma = \beta^{-2}$ in which case $a \in \beta^{-1}L^\perp(\gamma).$ 
Now let $L^\perp(\gamma) = \langle \alpha \rangle, \alpha \in A,$  and suppose we have $a = \beta^{-1}\alpha^k.$  
Then $\varphi(\beta \cdot \gamma) \sim \beta\cdot a\gamma=\alpha^k\gamma.$ 
Since $(\beta\gamma)^2=\beta^3$, squaring both sides gives $\beta^{3} \sim \beta$ which is a contradiction.
Thus we must have $a \in L^\perp( \gamma).  $
\qed 
 
\medskip


\begin{lem}\label{lem_phi(ar)} 
Let  $G$ be a wallpaper group with a reflection or glide reflection coset $Ar.$ 
Assume $\varphi|_A = \Id_A$ and $\varphi(Ar) = Ar.$  Suppose that for $ar \in Ar,$ we have $\varphi(ar) \sim ar.$  
Then $\varphi(ar) \in K_rar.$  
In particular, $ar = \varphi(r) \sim r$  implies that $a \in K_r.$  
\end{lem}
\bigskip

\noindent {\it Proof.} \indent    
If $\rho_\pi \notin G$, then by Lemma \ref {lemK}   we have $(ar)^G \cap Ar = K_rar$ and we are done, 
so suppose $\rho_\pi \in G.$  
Let $\lambda=(\rho_\pi, r^{-1})$, so that $\lambda \in A$ and 
$r^{\rho_\pi} = \lambda r.  $  
Let $\beta = r^2$ so that $r  = \beta r^{-1}.$
(If $r$ is a reflection then $\beta=1$.) 
Thus, again by Lemma \ref {lemK}   we have  \[(ar)^G \cap Ar  = K_rar \cup K_r(ar)^{\rho_\pi} = K_r ar \cup K_r a^{-1} \lambda r.\]
First we note that if $\lambda =1$ and $a^2 \in K_r$, then  we have nothing to prove; so suppose that this  is not the case.  
Since $\varphi(Ar) = Ar$ it suffices to show that $\varphi(ar) \notin K_r(ar)^{\rho_\pi}.$  
Now $K_r$ is a cyclic group: $K_r=\langle \alpha\rangle$, where $\alpha \in A$. Also in this case one sees that $K_r \le L^{\perp}(r)$. 
Suppose to the contrary that $\varphi(ar) = \alpha^i(ar)^{\rho_\pi}$  for some $\alpha^i \in K_r.$  
Then for any $b \in A$ we have 
\[
 \varphi(b \cdot ar) \sim b\alpha^ia^{-1}\lambda r.
\]
Squaring both sides we have 
\[
\varphi(bar)^2\sim \varphi((bar)^2)=    \varphi(ba\beta(ba)^r) = ba\beta(ba)^r \sim (b\alpha^ia^{-1}\lambda\beta)(b\alpha^ia^{-1}\lambda)^r.
\]
Reordering the elements this becomes
\begin{equation}\tag {5.1} \label{eqn_phi(af)} 
\beta aa^rbb^r \sim \beta \alpha^i(\alpha^i)^ra^{-1}(a^{-1})^rbb^r\lambda \lambda^r.
\end{equation}
Recall that $L(r)$ is the line fixed by the action of $r.$   
Note that $aa^r$ is in $L(r)$ (since $(aa^r)^r = a^ra = aa^r.$) 
Similarly $bb^r,  \lambda \lambda^r \in L(r).$  
Recall that  $\alpha \in K_r\le L^{\perp}(r), $ so that  $\alpha^r=\alpha^{-1}$.  So we know $ \alpha^i (\alpha^i)^r  =1.  $
 Also, $\beta \in L(r)$ (since $\beta^r = (r^2)^r = r^2 = \beta$).  
Thus Equation (\ref{eqn_phi(af)}) is stating that two elements of $L(r)$ are conjugate to each other.  This implies that they are equal or they are inverses.
If they are inverses of each other, Equation (\ref{eqn_phi(af)}) becomes
$b^{-2}(b^{-2})^r = \beta^2 \lambda \lambda^r $ which can't be true for all $b \in A.$ 
 If they are equal to each other, we have $a^2(a^2)^r = \lambda \lambda^r.$  
We can assume $\lambda$ is not trivial in this case because if $\lambda=1$,  then $a^2 \in K_r$, and we assumed both of those cannot be true together.  

So  assume $\lambda \ne 1.$ Now the fact that $G$ has a reflection or glide reflection $r$,  contains $\rho_\pi$, and $\lambda=(\rho_\pi,r^{-1})\ne 1$, restricts us to the groups of type {\bf p2mg, p2gg, p4mg}. In these cases one checks that
 $a^2(a^2)^r \in \langle x^4, y^4 \rangle$. However  (looking at the  possible   $\lambda$s in each of these three groups) we find that $\lambda \lambda^r \notin \langle x^4, y^4 \rangle $. This is a contradiction.
 Thus $\varphi(ar) \notin K_r(ar)^r.$    
\qed 
\medskip






\noindent   {\bf c2mm} case: Suppose that $\varphi(\rho)=x^uy^v\rho, \varphi(\sigma)=x^iy^j\sigma$. Then the fact that $(\rho\sigma)^2=1$ implies that
$1=\varphi((\rho\sigma)^2)=\varphi(\rho\sigma)^2=(\varphi(\rho)\varphi(\sigma))^2=(x^uy^v\rho x^iy^j\sigma)^2$, which is only true when $u+j=i+v$, so that we can compose $\varphi$ with the inverse of  $\psi_{u,v,i,j}$  to get $\varphi(\rho)=\rho,\varphi(\sigma)=\sigma.$ 
To get $\varphi(\rho\sigma)=\rho\sigma$ we now apply Lemma \ref {lemcfg}, and we are done.\medskip


\noindent   {\bf p2mm} case: Suppose that $\varphi(\rho)=x^uy^v\rho, \varphi(\sigma)=x^iy^j\sigma$. We can compose $\varphi$  with the inverse of  $\psi_{u,v}$  to get $\varphi(\rho)=\rho.$  To get $\varphi(\sigma) = \sigma$ and  $\varphi(\rho\sigma)=\rho\sigma$ we now apply Lemma \ref {lemcfg}, and we are done.\medskip



\noindent   {\bf p2mg} case: Suppose that $\varphi(\rho)=x^uy^v\rho.  $ Composing with a power of $\psi_y$ and then a power of $\psi_x$ we can arrange that $v=0$ and $u=0.$
Let $ \varphi(\sigma)=a\sigma$.
The fact that $\sigma^2 = 1$ gives $a \in L^\perp (\sigma).$  
In this group this implies $\varphi(\sigma) \sim \sigma$ and we may apply Lemma \ref{lem_phi(ar)} to give $a \in K_\sigma.$  
If $\varphi(\rho \sigma) = b\rho \sigma,$ then  Lemma \ref{lem_phi(gamma)} gives $b \in L^\perp (\rho \sigma)  $.  
Since 
$  (c\rho \sigma)^G = \langle x^2 \rangle ac \rho \sigma \cup \langle x^2 \rangle (cy)^{-1} \rho \sigma$,
$L^\perp (\rho \sigma)=\langle x^2\rangle$ and 
$b \rho \sigma = \varphi(\rho\cdot \sigma) \sim a^\rho \rho \sigma,$ we see  $a^\rho$ lies on 
$L^\perp (\rho \sigma) \cup y^{-1}L^\perp (\rho \sigma).$ Thus $a$ lies on $(L^\perp (\rho \sigma) \cup y L^\perp (\rho \sigma))\cap K_\sigma$. 
Hence $a =1,$ which gives us $\varphi(\rho \sigma) \sim \rho \sigma,$ and so by Lemma \ref{lem_phi(ar)} we know $b \in K_{\rho \sigma}.$  
Now $\sigma = \varphi(\rho \cdot \rho \sigma) \sim b^\rho\sigma$ which tells us 
$b \in L^\perp(\sigma), $ and so $b $ is also trivial.

\medskip 

\noindent   {\bf p2gg} case: Suppose that $\varphi(\rho)=x^uy^v\rho $. 
By using $\psi_x,\psi_y $ we can arrange that $u=v=0$.
Let $\varphi(\gamma) = a \gamma, \varphi(\rho \gamma) = b \rho \gamma.$  
By Lemma (\ref{lem_phi(gamma)}) $a \in L^\perp( \gamma)  .$  
Note that $L^\perp( \gamma)$ contains elements from only two conjugacy classes, namely $ \gamma^G $ and $ (y\gamma)^G .$  
If $a \gamma \sim y \gamma$ then we may compose with $\psi \circ \iota.$  
So now we   have $\varphi(\gamma) \sim \gamma$ and so by Lemma \ref{lem_phi(ar)} $a \in K_\gamma.  $
Lemma \ref{lem_phi(gamma)} gives us $b \in L^\perp(\rho \gamma)$ and thus $b \rho \gamma = \varphi(\rho \cdot \gamma) \sim            \rho a \gamma =a^{-1} \rho \gamma $ 
 implies that $a \in L^\perp(\rho \gamma) \cup yL^\perp(\rho \gamma).$  This intersects $K_\gamma $ trivially hence $\varphi(\gamma) = \gamma.$  

Now $\varphi(\rho \cdot \gamma) \sim \rho \gamma$ so we may apply Lemma \ref{lem_phi(ar)} which gives us $b \in K_{\rho\gamma}.$  
We also know $ \gamma = \varphi(\rho \cdot \rho \gamma) \sim b^{-1}\gamma$
which tells us $b \in K_\gamma \cup xyK_\gamma.$  
The only possibility then is that $b = 1.$

\medskip


\noindent   {\bf p3m1} case: Suppose that $\varphi(\rho)=x^uy^v\rho.$  
Composing with an element of  $\langle \psi_1, I_x\rangle$ we can arrange that $\varphi(\rho)=\rho$. 
Thus $\varphi(\rho)=\rho, \varphi(\rho^2)=\rho^2$. We now apply Lemma \ref {lemcfg} to conclude that $\varphi(\sigma)=\sigma, \varphi(\rho\sigma)=\rho\sigma, \varphi(\rho^2\sigma)=\rho^2\sigma$,
and we are done.
\medskip

\noindent   {\bf p31m} case: Suppose that $\varphi(\rho)=x^uy^v\rho, \varphi(\sigma)=x^iy^j\sigma$. Here we note that by composing with an element of  $\langle I_x, I_y\rangle$ we can arrange that $\varphi(\rho)=\rho$, so that $\varphi(\rho^2)=\rho^2$. From Lemma \ref {lemcfg} we obtain $\varphi(\rho^k\sigma)=\rho^k\sigma,k=0,1,2$, as required. 
\medskip

\noindent   {\bf p4mm} case: Suppose that $\varphi(\rho)=x^uy^v\rho, \varphi(\sigma)=x^iy^j\gamma$. 
Then acting by some $I_b, b \in A$, we can assume    $\varphi(\rho) \in \{\rho,y\rho\}$. If we have $\varphi(\rho)=y\rho$, then acting by the automorphism $\psi_1^{-1}$   we see that  $\varphi(\rho)=\rho, \varphi(\rho^3)=\rho^3$. Now Lemma \ref {lemcfg} gives $\varphi(\rho^k\sigma)=\rho^k\sigma, 0\le k\le 3$.
 
 Let $\varphi(\rho^2)=x^iy^j\rho^2$. Then for any reflection $t \in \{\sigma,\rho\sigma,\rho^2\sigma\rho^3\sigma\}$ we have
 $\rho^2t=\varphi(\rho^2t)\sim \varphi(\rho^2)\varphi(t)=x^iy^j\rho^2 t$, from which we see that $x^iy^j \in L^\perp (\rho^2 t)$ for all $t$. This gives $i=j=0$ and 
 so we   have $\varphi(\rho)=\rho, \varphi(\rho^2)=\rho^2, \varphi(\rho^3)=\rho^3, \varphi(\sigma)=\sigma$.
 By Lemma \ref {lemcfg} we   have $\varphi(\rho\sigma)=\rho\sigma, \varphi(\rho^2\sigma)=\rho^2\sigma, \varphi(\rho^3\sigma)=\rho^3\sigma.$
 
 So let $\varphi(\rho^2\sigma)=x^ay^b\rho^2\sigma$.  Then $\rho^{2+k}\sigma=\varphi(\rho^{2+k}\sigma)=\varphi(\rho^2\cdot \rho^k\sigma) \sim x^ay^b \rho^{2+k}\sigma$ gives $x^ay^b \in L^\perp(\rho^{2+k}\sigma)$ for $0\le k<4$. It follows that $a=b=0$, and this 
 concludes this case.\medskip

\noindent   {\bf p6m} case: Suppose that $\varphi(\rho)=x^uy^v\rho, \varphi(\sigma)=x^iy^j\gamma$. 
Then acting by some $I_b, b \in A$, we can assume that   $\varphi(\rho) =\rho$.  
Now Lemma \ref {lemcfg} gives $\varphi(\rho^k\sigma)=\rho^k\sigma, 0\le k\le 5$, and the rest of the proof follows as in the previous case.
 This completes this case and concludes the proof of Step (4). 

\section {Steps (5) and (6)} 

\begin{lem} \label {lemra=rb} Suppose that $\varphi|_A=\Id_A$ and $t \in G$ where $\varphi(t)=t$. If $a \in A$ where $\varphi(at)=bt$, then $a \sim b$.
\end{lem}
\noindent {\it Proof} 
Given the hypotheses we have:
$$a \sim a^t=\varphi(a^t)=\varphi(t^{-1}\cdot at)\sim \varphi(t^{-1})\cdot \varphi(at)=t^{-1}bt\sim b.\qed $$

Recall that $F$ is a set of coset representatives for $G/A$.
In the situation of Lemma \ref {lemra=rb} we will write $\varphi(at)=a^{r_a}t$, where $r_a \in F$.

\begin{prop} \label{propra=rb}
Suppose that $\varphi|_A=\Id_A$ and $t \in G$ where $\varphi(t)=t$. Then there is  $f \in F$ such that $\varphi(at)=a^ft$ for all $a \in A$. 
\end{prop}
\noindent {\it Proof} From Lemma \ref {lem2.0} we see that $\varphi(At)=At$. 
 Now for $a,b \in A$ we have
 $$ab^{-1}=at
 (bt)^{-1}=\varphi(  at\cdot (bt)^{-1} ) \sim a^{r_a}t\cdot t^{-1}(b^{r_b})^{-1}=a^{r_a}(b^{r_b})^{-1}.
 $$ Thus there is some $f \in F$ such that $ab^{-1}=(a^{r_a}(b^{r_b})^{-1})^f$, so that letting $\alpha=r_af, \beta=r_bf$ we have
\begin{align*}\tag {6.1} \label {tag6.1}
a (a^{-1})^\alpha  = b (b^{-1})^\beta.
\end{align*}

For $a \in A$ we let $C_a$ denote the circle in $\mathbb E^2$ that contains the  origin and  that is centered at $v_a$, and let $S_a=\{v_{a (a^{-1})^f}:f \in F\}$. Then $S_a$ consists of  $|a^G|$ points that  lie on $C_a$.  Note that $(0,0) \in S_a$ for all $a \in A$.  Since two distinct circles can meet in at most  two points  equation ({\ref {tag6.1})  gives
$|S_a \cap S_b| \in \{1,2\}$.

\begin{lem} \label {lemkl}
(i) If $|S_a \cap S_b|=1$, then there is $\delta \in F$ such that $a^{r_a}=a^\delta, b^{r_b}=b^\delta$. 

(ii) If $(0,0),v_a,v_b$ are collinear, then there is a $\delta \in F$ such that $a^{r_a}=a^\delta, b^{r_b}=b^\delta$. 
\end{lem}
\noindent {\it Proof}
(i) If $|S_a \cap S_b|=1$, then  $S_a \cap S_b=\{(0,0)\}$ and so from equation (\ref {tag6.1}) we get $a (a^{-1})^\alpha  = b (b^{-1})^\beta=1$, so that $a =a^\alpha,  b =b^\beta.$ From this we obtain 
$a^{f^{-1}}=a^{r_a}, b^{f^{-1}}=b^{r_b}$, so that we can let $\delta=f^{-1}$.

 (ii) We may  assume that $a \ne b$.
Since $(0,0),v_a,v_b$ are collinear the centers of the circles $C_a, C_b$ are on the line through these points. Since $a \ne b$ we have $C_a \ne C_b$; since $(0,0)$ is common to $C_a$ and $C_b$ we see that $C_a \ne C_b$ implies  $C_a\cap C_b =\{(0,0)\}$. This gives $a =a^\alpha,  b =b^\beta,$ which then implies 
$a^{f^{-1}}=a^{r_a}, b^{f^{-1}}=b^{r_b}$, as required.\qed \medskip 

\begin{lem} \label {lembnotinH} Let $b \notin H$. Then for all $a \in A$ there is a unique $f=f_{a,b}  \in F$ such that $\varphi(at)=a^ft, \varphi(bt)=b^ft$.
\end{lem}
\noindent {\it Proof} 
Recall that $b \notin H$  means that $|b^G|=|F|$.  If  $(0,0),v_a, v_b$ are collinear, then the existence of such a $\delta\in F$ 
 follows from Lemma \ref {lemkl} (ii), while the uniqueness follows from the fact  that $b \notin H$.
 
Now assume   that $(0,0),v_a, v_b$ are not collinear.  If  $|S_a \cap S_b|=1$, then by Lemma \ref {lemkl} (i) there is $\delta \in F$ such that $a^{r_a}=a^\delta, b^{r_b}=b^\delta$, and the fact that $\delta$ is unique follows from $b \notin H$. 
So now assume that $|S_a \cap S_b|=2$. 
Then by Lemma \ref {lemkl} (ii) there is some $f \in F$ such that for all $k\in \mathbb N$ we have
$(b^k)^f=(b^k)^{r_{b^k}}$. Since $b \notin H$ this element $f$ is unique.
Since $S_a$ is finite there is some $k \in \mathbb N$ such that $|S_a \cap S_{b^k} |=1$; then we have $h,h' \in F$ with $b^{r_b}=b^h, b^{r_{b^k}}=(b^k)^h$
and $ b^{r_{b^k}}=(b^k)^{h'}, a^{r_a}=a^{h'}.$ Since $b \notin H$ we again see that $h,h'$ are unique, so that $h=h'$. Thus we have 
$ a^{r_a}=a^{h}$ and $b^{r_b}=b^{h}$, as required. \qed \medskip

Now let $a,b,c \in A$ where $b \notin H$. Then by Lemma \ref {lembnotinH} there are unique $f,h \in F$ such that $\varphi(at)=a^ft, \varphi(bt)=b^ft$ and 
$\varphi(bt)=b^ht, \varphi(ct)=c^ht$. Since $f,h $ are unique and $ \varphi(bt)=b^ft,  \varphi(bt)=b^ht$ we must have $f=h$; it follows (by fixing $a$ and varying $c$) that for all $d \in A$ we must have $\varphi(dt)=d^ft$ for this value of $f \in F$ that is completely determined by $b \notin H$.  
This concludes the proof of Proposition \ref {propra=rb}, and so completes Step (5).\qed\medskip

For the remainder of this section we assume that  $\varphi|_A=\Id_A$ and that $\varphi(f)=f$ for all $f \in F$. We also assume that for all $f \in F$ there is $r_f$ such that $\varphi(af)=a^{r_f}f$ for all $a \in A$.  Our goal will be to reduce  to the situation where 
 $\varphi(af)=af$ for all $f \in F,a \in A$.  We start with:



\begin{lem}\label{lem_1_or_r}
 For all $a \in A$ and all reflections or glide reflections $r \in F$ we have $\varphi(ar)\in \{ar, a^r r\}$.
\end{lem}

\noindent {\it Proof.} \indent
Since  $ar \sim \varphi(a\cdot r)$  we may apply Lemma  \ref{lem_phi(ar)}.
This gives us $\varphi(ar) = a^{t}r \in K_rar$ for some $t \in F.$  This implies that $a^{-1}a^{t} \in K_r$ for all $a \in A.$ 
But $a^{-1}a^t=(a,t)$ and $K_t=\{(a,t):a \in A\}$ shows that $K_t \subset K_r$. From Lemma \ref {lemK} (iii) it follows that either $K_t=\{1\}$ or $K_t=K_r$, from which we see that either $t=1$ or $t=r$.
\qed
\medskip

The following applies to the   groups of type  {\bf  p4mm, p4mg, p31m, p3m1,  p6m}.
\begin{thm}\label{thmIdonRefl}
 Let $G$ be a wallpaper group having a nontrivial rotation $\rho\neq \rho_\pi$, and distinct reflection or glide reflection cosets $Ar, Ar \rho$.  Suppose we have 
 $\varphi |_A=\Id_A,$ 
 $\varphi(r)=r, \varphi(r\rho) =r \rho,$ and $\varphi(\rho) = \rho.$ 
Then $\varphi|_{Ar}= \Id_{Ar}.$   
 \end{thm}
\noindent {\it Proof.} 
Recall from Lemma  \ref {lemK}    that for any of these five groups we have
\[
(ar)^G \cap Ar \subseteq K_r ar \cup K_r (a\beta_1 \beta_2)^{-1}r,
\]
where $\beta_1\beta_2 = xy$ when $G$ is of type {\bf p4mg}, but is trivial otherwise. 
Suppose to the contrary that $\varphi|_{Ar}\ne \Id_{Ar}$. Then by Lemma \ref {lem_1_or_r} we must have $\varphi(ar)= a^rr$ for all $ a \in A$. Now 
$r \rho$ is also a reflection or a glide reflection, and so we have:
\[
a^r r\cdot\rho\sim \varphi(ar \cdot \rho) =\varphi(a \cdot r\rho) \sim ar\rho. 
\]
Thus for all $a \in A$ we must have 
(i) $a^r a^{-1} \in K_{r \rho}$, or (ii) $a^r a \in K_{r \rho}(\beta_1\beta_2)^{-1}.$
Since $r$ and $r\rho$ are distinct reflections, case (i) is only true for $a$ which satisfy $a^r a^{-1} = 1.$  Any such $a \in A$ must lie on $L(r).$ 
Since $\rho\ne \rho_\pi$, $L(r) $ and $K_{r\rho}$  are not  parallel.  Hence $L(r) \cap K_{r\rho}(\beta_1\beta_2)^{-1}$ is exactly one point. 
(If $G$ is of type {\bf p4mg} and  $r = \rho \gamma$ this point is $xy^{-1}$, when $r = \rho^2 \gamma$ it is $y^{-2},$ and when $r= \rho^3 \gamma$ it is $x^{-1}y^{-1}.$  In all other cases, this point is trivial.)
Thus case (ii) is only true for $a \in L^\perp(r)y^{-1}$ or  $a \in L^\perp(r)$.   
This shows that the only $a \in A$ that satisfy $a^rr\rho \sim ar\rho$ are in the union of two lines (either $L(r) \cup L^\perp(r)$ or $L(r) \cup L^\perp(r)y^{-1}$) which is not all of $A.$  
Thus $\varphi(ar) = a^r r $ is not possible.
We conclude that $\varphi|_{A r} = Id.$ 
\qed\medskip
   
 The following applies to the eight groups {\bf c2mm, p2mm, p2mg, p2gg, p4mm, p4mg, p31m, p3m1, and p6mm}.


\begin{thm}\label{thmIdonRot}
 Let $G$ be a wallpaper group having at least two cosets that are reflection or glide reflection cosets. Suppose we have 
 $\varphi |_A=\Id_A$ and $\varphi|_{A\sigma} = Id$ for all reflection cosets $A \sigma$ as well as $\varphi|_{A\gamma} = Id$ for all glide reflection cosets $A \gamma.$   
Assume $\varphi(t)=t$ for all $t \in F$.
Then $\varphi= Id$.
 \end{thm}
\noindent{\it Proof}
Let $Ar$ denote a reflection or glide reflection coset.  
Let $A\rho^k$ be a rotation coset.  We know by Proposition \ref {propra=rb} that there is some $t \in F$ such that 
$\varphi(a\rho^k)=a^t\rho^k$ for all $a \in A$. By hypothesis we have $\varphi(r a^{-1})=r a^{-1}$ which gives
\[
r\cdot \rho^k\sim \varphi(r \cdot \rho^k) = \varphi(r a^{-1}\cdot a\rho^k)\sim r a^{-1}\cdot a^t\rho^k=(a^{-1}a^t)^r \cdot r\rho^k.
\]
(The last equality is true because the action of $r$ is the same as the action of $r^{-1}.$)
This gives $(a^{-1}a^t)^r  \in K_{r\rho^k} \cup K_{r\rho^k}(\beta_1\beta_2)^{-1}$ for all $a \in A.$ Write $\beta_r = ((\beta_1\beta_2)^{-1})^{r^{-1}}$ so that $(a,t) \in  K_{\rho^kr}\cup K_{\rho^kr}\beta_r $ for all $a \in A$. 
Now the element $t \in F$ depends upon $k$ and not on $r$. Then the fact that we have at least two reflection or glide reflection cosets, say $r, r^\prime \in F$ means that
  $(a,t)\in \big( K_{\rho^kr} \cup K_{\rho^kr}\beta_r\big)\cap  \big( K_{\rho^kr^\prime} \cup K_{\rho^kr^\prime}\beta_{r^\prime}\big).$
This intersection is trivial if neither $r$ nor $r^\prime$ are glide reflections.  
Otherwise, this intersection is finite and therefore bounded. 
Thus there exists a $t \in F$ such that for all $a \in A,$ we must have $ (a,t) $   in this small intersection. 
This can only be true if $t$ is trivial.  
So   $\varphi(a\rho^k)=a\rho^k$. This completes the proof.\qed\medskip

\noindent {\it The group of type} {\bf c2mm}  
By Lemma \ref{lem_1_or_r} there are $s \in \{1, \sigma \}, u \in \{1, \rho \sigma \}$, where
$ \varphi(a\sigma)=a^s\sigma, \varphi(a\rho\sigma)=a^u\rho\sigma$ for all $a \in A.$  
If we have $u = \rho \sigma,$ we may compose with $\psi \circ \iota$ so that we have $u=1.$ If  $s = \sigma$, then  we have
\[
x y \rho \sim \varphi(xy \cdot \rho) = \varphi(x \sigma \cdot x \rho \sigma) \sim    y \sigma x \rho \sigma = y^2 \rho .
\]
Since $xy \rho \nsim \,  y^2 \rho,$  this is not possible so we must have  $s = 1.$  Theorem \ref {thmIdonRot} now  shows that $\varphi(a\rho)=a\rho$ for all $a \in A$, and so $G$ has trivial weak Cayley table group.\qed 
 
\medskip


\noindent {\it The group of type} {\bf p2mm} 

By Lemma \ref{lem_1_or_r} there are $s \in \{1, \sigma \}, u \in \{1, \rho \sigma \}$, where
$ \varphi(a\sigma)=a^s\sigma, \varphi(a\rho\sigma)=a^u\rho\sigma$ for all $a \in A.$  
If we have $s = \sigma,$ we may compose with $\tau$ so that we have $s=1$.

Now note that $\tau^\psi$ is the identity on $A \cup A\sigma$ and conjugates elements of $A\rho \cup A\rho\sigma$ by $\rho\sigma$.
This, 
if we now have $u=\rho\sigma$, then we compose with $\tau^\psi$ to obtain $s=1,u=1$. Theorem \ref {thmIdonRot} now  shows that $\varphi(a\rho)=a\rho$  for all $a \in A$. This shows that groups $G$ of type  {\bf p2mm} have weak Cayley table group generated by $\mathcal W_0(G)$ and  $\tau$. \qed

\medskip 

\noindent {\it The group of type } {\bf p2mg}  
By Lemma \ref{lem_1_or_r} we know that 
$ \varphi(a\sigma)=a^s\sigma,$ and $ \varphi(a\rho\sigma)=a^u\rho\sigma$ for some $s \in \{1, \sigma \}, u \in \{1, \rho \sigma \}$, for all $a \in A.$
Now if $s=\sigma$, then
$$\rho\sigma=\varphi(\rho\sigma)=\varphi(y\sigma\cdot \rho)\sim y^s\sigma\rho=y^sy^{-1}\rho\sigma=y^{-2}\rho\sigma,
$$ which contradicts Lemma \ref {lemK}. Thus $s=1$.

Next, if $u=\rho\sigma$, then
$$x\sigma=\varphi(x\sigma)=\varphi(x\sigma\rho\cdot\rho)=\varphi(xy^{-1}\rho\sigma\cdot \rho)\sim (xy^{-1})^u\rho\sigma\rho=
 (xy^{-1})^uy\sigma,$$
 so that $ (xy^{-1})^uy=x^{-1}y^{-1}y=x^{-1}$ has the form $xy^{2k}$ or $x^{-1}y^{2k}y^{-1}$, both of which are not possible. Thus $u=1$.
  Theorem \ref {thmIdonRot} now  shows that $\varphi(a\rho)=a\rho$ for all $a \in A$, and so $G$ has trivial weak Cayley table group.\qed\medskip

\noindent {\it The group} {\bf p2gg}  
By Proposition  \ref {propra=rb}  there are $s,u \in F=\{1,\rho,\gamma,\rho\gamma\}$ such that
$\varphi(a\gamma)=a^s\gamma,$ and $ \varphi(a\rho\gamma)=a^u\rho\gamma$ for all $a \in A$.  
By Lemma \ref{lem_1_or_r} we know $s \in \{1, \gamma\}, u \in \{1, \rho \gamma\}.$  
Let $X=\langle x^2\rangle, Y=\langle y^2\rangle$. Then for $G$ a group of type  {\bf p2gg}  we have

\begin{lem}\label {lemclp2gg} Let $a,b \in A$. Then

\noindent (i) $a\gamma \sim b \gamma$ if and only if $ab^{-1} \in Y $ or $ab \in x^{-1}yY$.

\noindent (ii) $a\rho\gamma \sim b \rho\gamma$ if and only if $ab^{-1} \in X$ or $ab \in xy^{-1}X.$

\noindent (iii) If $v \in F$ and $a \in A$, then $aa^v \notin x^{-1}yY \cup xy^{-1}X$.

\end{lem}
\noindent {\it Proof} (i) and (ii) follow from the presentation. 
For (iii) note that for $a=x^iy^j$ we have $a^\rho=a^{-1}, a^\gamma=x^iy^{-j}, a^{\rho\gamma}=x^{-i}y^j$, which gives (iii).\qed

Now for all $a \in A$ we have:
\begin{align*} (a^{-1})^{\gamma u}\rho\gamma&=
a^{\rho\gamma u}\rho\gamma=\varphi(a^{\rho\gamma}\rho\gamma)=\varphi((a\rho\gamma)^{\rho\gamma})\sim \varphi(a\rho\gamma)\sim 
 \varphi(\rho\gamma\cdot a )\\
 &=\varphi(\rho a^{\gamma^{-1}}\gamma)=\varphi(\rho \cdot a^{\gamma}\gamma)
 \sim \rho (a^\gamma)^s\gamma = a^{\gamma s\rho}\rho\gamma=(a^{-1})^{\gamma s}\rho\gamma.
\end{align*}

From Lemma \ref {lemclp2gg} (ii), (iii), we thus get $a^s(a^{-1})^u \in X$ for all $a \in A$.
Let $a=x^iy^j$, and recall that $u \in \{1,\rho\gamma\}$.
We show that this implies  $s=1$, for if $s=\gamma$, then the possibility $u=1$ gives $a^s(a^{-1})^u=(x^iy^{-j})(x^{-i}y^{-j})=y^{-2j}\notin X$; while  $u=\rho\gamma$ gives $a^s(a^{-1})^u=(x^iy^{-j})(x^iy^{-j})=x^{2i}y^{-2j}\notin X$. Thus we  have $s=1$. 

Similarly, for all $a \in A$ we have:
$$a^{\gamma }\gamma=\varphi(a^\gamma \gamma)=\varphi(\rho\cdot a^{\rho \gamma}\rho\gamma)\sim \rho a^{\rho\gamma u}\rho\gamma
=a^{\gamma u}\gamma,
$$
from which we obtain $(a^{-1})^\gamma a^{\gamma u} \in Y$ for all $a \in A$; or  $a^{-1} a^{ u} \in Y$ for all $a \in A$. However, if $u=\rho\gamma$, then $x^{-1}x^{\rho\gamma}=x^{-2}$, a contradiction. Thus $u=1$.  
\qed\medskip

  Theorem \ref {thmIdonRot} now  shows that $\varphi(a\rho)=a\rho$ for all $a \in A$, and so any $G$ of type {\bf p2gg} has trivial weak Cayley table group.\qed\medskip

Using  the above cases,  Theorem \ref {thmIdonRefl} and Theorem  \ref {thmIdonRot}  we see that 

\begin{cor} \label {corngl} Each group of type  {\bf   c2mm, p2mg, p2gg, p4mm, p31m, p3m1, p6mm} has trivial weak Cayley table group.\qed\end{cor} 

\noindent {\it The group} {\bf p2} 

We have  $\varphi|_A=\Id_A$ and that there is some $t \in F$ such that $\varphi(a\rho)=a^t\rho$ for all $a \in A$. If $t=1$, then $\varphi=id$, while if $t=\rho$, then $\varphi= I_\rho\iota$, and we have proved

\begin{cor} \label {cornglp2} The group of type  {\bf p2} has trivial weak Cayley table group.\qed
\end{cor} 

\medskip 

\begin{lem} \label {lemrrr}  Let $G$ be one of the groups of type {\bf p3, p4, p6} (respectively), and let $d=3,4,6$ (respectively). 
Suppose that $\rho \in F$ is a rotation of order $d$, and that there are $s,t \in F$ with  $\varphi(a\rho)=a^t\rho, \varphi(ar^{-1})=a^s\rho^{-1}$ for all $a \in A$. Then $s=t$.\end{lem}
\noindent {\it Proof}  
Given $a \in A$,  let $a\rho^{d-1}=(b\rho)^{-1}, b \in A$. Then solving gives $b=(a^{-1})^{\rho^2}$. Since $(t,\rho)=1$ we have:
$$\varphi(a\rho^{d-1})=\varphi(b\rho)^{-1}=(b^t\rho)^{-1}=\rho^{-1}(b^{-1})^t
=\rho^{-1}(a^{\rho^{-1}})^t=a^t\rho^{d-1}\qed \medskip
$$

\noindent {\it The group} {\bf p3} 

We know that $\varphi|_A=\Id_A$ and that there is  $t \in F=\{1,\rho,\rho^2\}$ such that $\varphi(a\rho)=a^t\rho$ for all $a \in A$. By Lemma \ref {lemrrr} we have $\varphi(a\rho^2)=a^t\rho^2$ for all $a \in A$. 
 Thus $\varphi$ is a power of the weak Cayley table isomorphism $\tau$ defined in $\S 2$, and we have proved

\begin{cor} \label {cornglp3} The group of type  {\bf p3} has   weak Cayley table group generated by $\mathcal W_0(G)$ and $\tau$.\qed
\end{cor} 

\medskip 

\noindent {\it The group} {\bf p4} 

We know that $\varphi|_A=\Id_A$ and that there is some $t \in F=\langle \rho\rangle$ such that $\varphi(a\rho)=a^t\rho$ for all $a \in A$. By Lemma \ref {lemrrr} we have $\varphi(a\rho^3)=a^t\rho^3$ for all $a \in A$. If $t=1$, then $\varphi|_{A\rho\cup A\rho^3}=id.$ If $t=\rho^k$, then composing with   $\mu_{\rho}^{-k}$ reduces to the case where 
 $\varphi|_{A\rho\cup A\rho^3}=id$, also. 

Now there is $u \in F$ such that $\varphi(a\rho^2)=a^u\rho^2$ for all $a \in A$. If $u =\rho^2$, then we can compose with $\tau_{\rho^2}$ to obtain $\varphi=id$.
Thus the remaining case reduces to $u=\rho$. In this situation one has $\varphi(x^2\rho^2)=(x^2)^\rho\rho^2=y^2\rho^2$, while $\varphi(x\cdot x\rho^2)\sim x\cdot x^\rho \rho^2=x\cdot y\rho^2$, giving $y^2\rho^2 \sim xy\rho^2$, a contradiction. This gives

\begin{cor} \label {cornglp4} The group of type  {\bf p4} has   weak Cayley table group generated by $\mathcal W_0(G)$ and $\tau_x, \tau_y, \tau_{\rho^2}, \mu_{\rho}$.\qed
\end{cor}

\noindent {\it The group} {\bf p6}

We know that $\varphi|_A=\Id_A$ and that there is some $t \in F=\langle \rho\rangle$ such that $\varphi(a\rho)=a^t\rho$ for all $a \in A$. By Lemma \ref {lemrrr} we have $\varphi(a\rho^5)=a^t\rho^5$ for all $a \in A$. Similarly, there is $s \in F$ such that if  $\varphi(a\rho^2)=a^s\rho^2$ for all $a \in A$, then we also have 
 $\varphi(a\rho^4)=a^s\rho^4$ for all $a \in A$.
 Lastly, assume that  $\varphi(a\rho^3)=a^u\rho^3$ for all $a \in A$.
 
 Now using an element of $\langle \tau_{\rho^2}I_{\rho^2},\mu_{\rho^3}I_{\rho^3}\rangle$  we can reduce to the case where 
 $t=1$ while preserving the conditions $\varphi|_A=\Id_A, \varphi(\rho^k)=\rho^k, 0\le k<6$.

 \begin{lem} \label{lemc6hex}  Assume that for some $t \in \langle \rho\rangle$ we have $a^{-1}a^t \in K_{\rho^2} \cap K_{\rho^3}$ for all $a \in A$. Then $t =1$.
  \end{lem}
\noindent {\it Proof} If  $a=y$ then $\{a^{-1}a^v: v \in \langle \rho \rangle\}=\{1,x^{-1}, x^{-1}  y^{-1}, x  y^{-2}, x  y^{-1}, y^{-2}\}$. 
Now $$K_{\rho^2}=\langle (x,\rho^2),(y,\rho^2)\rangle=\langle x^{-2}y,xy\rangle\text { and } 
K_{\rho^3}=\langle x^2,y^2\rangle,$$ so that  $K_{\rho^2}\cap K_{\rho^3}=\langle x^2y^2, x^2y^{-4}  \rangle$.  But the only element of 
$\{1,x^{-1}, x^{-1}  y^{-1}, x  y^{-2},$ $ x  y^{-1}, y^{-2}\}$ that is also in $K_{\rho^2}\cap K_{\rho^3}$ is the identity, and so we are done.\qed\medskip

Now 

\noindent (i) $a^s\rho^2=\varphi(a\rho^2)=\varphi(a\rho\cdot \rho)\sim a\rho^2$ gives  $a^{-1}a^s \in K_{\rho^2}$; 

\noindent (ii) $a^s\rho^4=\varphi(a\rho^3\cdot \rho)\sim a^u\rho^4$ gives  $(a^{-1})^sa^u \in K_{\rho^2}$;

\noindent (iii) $a^u\rho^3=\varphi(a\rho^3)=\varphi(a\rho\cdot \rho^2)\sim a\rho^3$ gives $a^{-1}a^u \in K_{\rho^3}$;

\noindent (iv) $a^u\rho^3=\varphi(a\rho^2\cdot \rho)\sim a^s\rho^3$ give $a^u(a^{-1})^s\in K_{\rho^3}$.

Since $K_{\rho^2}, K_{\rho^3}$ are subgroups (i) and (ii) give $a^{-1}a^u \in K_{\rho^2}$ for all $a \in A$. Thus by (iii) we have
$a^{-1}a^u \in K_{\rho^2} \cap K_{\rho^3}$, so that Lemma \ref {lemc6hex} tells us that $u=1$.
 
Similarly, (iii) and (iv) give $a^{-1}a^s \in K_{\rho^3}$, which combined with (i) gives $ a^{-1}a^s \in K_{\rho^2} \cap K_{\rho^3}$, from which we obtain $s=1$, using Lemma \ref {lemc6hex}. Thus  we now have $\varphi=\Id_G$. This gives:

\begin{cor} \label {cornglp6} The group of type  {\bf p6} has   weak Cayley table group generated by $\mathcal W_0(G)$ and $\tau_{\rho^2}, \mu_{\rho^3},\mu_{x^2},\mu_{y^2}$.\qed
\end{cor} 

\medskip

\noindent {\it The groups} {\bf cm, pm, pg}  

Here we have $t \in F$ such that $\varphi(a\sigma)=a^t\sigma$. If $t=1$, then we are done. So assume that $t\ne 1, t \in F$.
In each such group we have the automorphism  $\psi: (x,y,r) \mapsto (x^{-1},y^{-1},r^{-1})$, where $r$ is $\sigma, \sigma,\gamma$ respectively, for the three groups.
 Then composing $\psi$ with $\iota$ gives $\varphi$, and we have proved:

 \begin{cor} \label {cornglp6} The groups of type  {\bf cm, pm, pg} have trivial   weak Cayley table groups.\qed
\end{cor} 
\medskip 

\section{The groups of type {\bf p4mg}}


We note that much of what we have done in other cases does not apply to a group of type {\bf p4mg}.

Here we start by assuming  $\varphi|_A=\Id_A$ and  $\varphi(At)=At$ for all $t \in F$. 
Let $$U=\langle xy,x^2\rangle, \,\,\, X=\langle x^2\rangle, \,\,\,Y=\langle y^2\rangle, \,\,\,V=\langle xy\rangle, \,\,\,W=\langle xy^{-1}\rangle.$$
In the next two results we list the  conjugacy classes and involutions:

\begin{lem} \label {lemclp4mg} 
 Let $a=x^iy^j \in A$. Then 
 \begin{align*}
 \:& & (x^iy^j)^G &= \{ x^{\pm i}y^{\pm j} , x^{\pm i}y^{\pm j} \};\\
& & (a\rho)^G &= U a\rho \cup U ax\rho^3 ;\\
& & (a\rho^2)^G &= U a\rho^2 ;\\
& & (a \rho^3)^G & = U a \rho^3 \cup U ax \rho ;\\
& & (a\gamma )^G &= Y a\gamma \cup Y a^{-1}x^{-1}y\gamma \cup X a^\rho x\rho^2\gamma \cup X a^{\rho^3}y^{-1}\rho^2\gamma ;\\
& & (a\rho \gamma)^G &= V a\rho \gamma \cup V a^{-1}\rho \gamma \cup W y^{-1}a^\rho\rho^3\gamma \cup W x^{-1}(a^{-1})^\rho\rho^3\gamma ;\\
& & (a \rho^2 \gamma)^G &= X a \rho^2 \gamma \cup X a^{-1}xy^{-1} \rho^2 \gamma \cup Y a^{\rho^3} y \gamma \cup Y a^\rho x^{-1} \gamma ;\\
& & (a\rho^3 \gamma)^G &= W a\rho^3 \gamma \cup W x^{-1}y^{-1}a^{-1}\rho^3 \gamma \cup V xa^{\rho^3} \rho \gamma \cup V x^{-1}a^\rho \rho \gamma.\qed
\end{align*}\end{lem}
 
 \begin{lem} \label {leminvp4mg} Any involution in $G$ is in one of the cosets:
 $$(a) \,\, A\rho^2;\quad (b) \,\, V\rho\gamma; \quad (c) \,\, x^{-1}W\rho^3\gamma.\qed $$
 \end {lem}


Now assume that $\varphi(\rho)=a\rho, a \in A$. Note that $I_x(\rho)=(xy)^{-1}\rho, I_y(\rho)=xy^{-1}\rho$, so that by acting by some $I_b, b \in A$, we can assume    $\varphi(\rho) \in \{\rho,x\rho\}$. If we have $\varphi(\rho)=x\rho$, then composing with the automorphism $\psi_1$   we see that  $\varphi(\rho)=\rho$.  
Now let $\varphi(\gamma)=a\gamma$.  
 Lemma \ref{lem_phi(gamma)} gives $a \in \langle y \rangle.$ 
 Note that $(\rho\gamma)^2=1$ and Lemma \ref   {leminvp4mg}  implies that $\varphi(\rho \gamma) \sim \rho \gamma.$ 
So Lemma \ref{leminvp4mg} (b) with $ \varphi(\rho \cdot \gamma) \sim \rho a \gamma 
$ tells us $a^{\rho^3} \in V$, so $a \in V^\rho=W.$  
We conclude that $a \in \langle y\rangle \cap W$ is trivial. 
Thus we now have $\varphi(\rho)=\rho, \varphi(\gamma)=\gamma$, so that  $\varphi(\rho^3)=\rho^3$.

Now let $\varphi(\rho^2\gamma)= b\rho^2\gamma, b \in A$. 
 Lemma \ref{lem_phi(gamma)} gives $b \in \langle x \rangle.$
Then using Lemma \ref{leminvp4mg} (b) again, $(\rho\gamma)^2=1$ and $\varphi(\rho\gamma)=\varphi(\rho^{-1}\cdot \rho^2\gamma)
\sim \rho^{-1}b\rho^2\gamma=b^\rho\rho\gamma$ 
gives $b^\rho \in V,$ hence $b \in W \cap \langle x\rangle.$  Thus $b=1$ and  $\varphi(\rho^2 \gamma) = \rho^2 \gamma.$  

Assume that $\varphi(\rho^3\gamma)=c\rho^3\gamma, c \in A$. 
By Lemma \ref{lem_phi(gamma)} we have $c \in W.$
Note that
$
\gamma = \varphi(\rho \cdot \rho^3 \gamma) \sim\rho c\rho^3\gamma= c^{\rho^3} \gamma$ implies that
$ c^{\rho^3} \in Y \cup Yx^{-1}y$, and so   $ c \in X \cup Xxy^{-1}.
$
But also,
$
\rho^2 \gamma = \varphi(\rho^{3} \cdot \rho^3 \gamma) \sim c^{\rho} \rho^2\gamma$ implies that  $ c^{\rho} \in X \cup Xxy^{-1},$ so that
 $ c \in Y \cup Yx^{-1}y.
$
As $c \in W \cap (X \cup Xxy^{-1}) \cap (Y \cup Yx^{-1}y)$  we see that $c=1$, and so $\varphi(\rho^3\gamma) = \rho^3 \gamma.$

The above shows that $\varphi(v)=v$ for all $v \in F \setminus \{\rho^2, \rho \gamma\}$. Thus by Proposition \ref {propra=rb}, for  $ v \in F \setminus \{\rho^2, \rho \gamma\}$ there is $t_v \in F \setminus \{ \rho^2, \rho \gamma \} $ such that  for all $a \in A$ we have
$\varphi(av)=a^{t_v}v$.
Let $t, u, w \in F$ satisfy
$\varphi(a\rho) = a^t \rho,
 \varphi(a\gamma) = a^u\gamma ,$ and 
$\varphi(a \rho^{-1} \gamma) = a^w \rho^{-1} \gamma.$
By Lemma \ref{lem_1_or_r}} we have $u \in \{1, \gamma \}$ and $w \in \{1, \rho^{-1}\gamma \}$.

Now let $\varphi(\rho\gamma)=a\rho\gamma, a \in A$. Then Lemma \ref {leminvp4mg} (b)  gives $a\in V$, while  $\varphi(\rho \cdot \rho\gamma) \sim a^{\rho^{-1}}\rho^2 \gamma$ tells us $a \in Y \cup Yxy$.  Therefore we know $a \in \{1, xy\}$.  
However if we write $t^{\rho \gamma} = s, $ and assume $\varphi(\rho \gamma) = xy\rho \gamma,$ then 
\[
x^2 \gamma \sim \varphi(x^3 \cdot x^{-1} \gamma) =
\varphi(x^3 \rho \gamma \rho ) =
\varphi(\rho \gamma \cdot (x^3)^{\rho \gamma} \rho) 
\]
\[
\sim xy \rho \gamma (x^3)^{\rho \gamma t} \rho
= xy (x^3)^{s} \rho \gamma \rho
= xy (x^3)^s \gamma^{-1}
= y (x^3)^s \gamma.
\]
This gives a contradiction unless $s \in \{\rho^2, \rho^2 \gamma \}$ which corresponds to $t \in \{\rho^2, \gamma \}.$  
But since $(y^{-1})^{\rho^2} = (y^{-1})^\gamma = y$ and 
$(x^{-1})^\gamma = x^{-1}$  for any combination of $t \in \{\rho^2, \gamma \}$ and $u \in \{ 1, \gamma \}$ we have
\[
\rho \gamma \sim \varphi(\rho \gamma) = \varphi(y^{-1}\rho \cdot x^{-1}\gamma) \sim (y^{-1})^t \rho (x^{-1})^{u}  \gamma 
=y \rho(x^{-1})    \gamma =y^2\rho \gamma,
\]
 which is a contradiction.  
Thus  $\varphi(\rho \gamma) = \rho \gamma.$  

Now we have $\varphi(\rho^k \gamma) = \rho^k \gamma$ for  $k \in \{0, 1, 2, 3\},$ and we may apply Theorem \ref {thmIdonRefl} four times to get $\varphi = Id$ on both reflection cosets and both glide reflection cosets.  

Let $\varphi(\rho^2) = b\rho^2, b \in A.$  
Since $\varphi(\rho^2 \cdot \gamma)  \sim b \rho^2 \gamma$ we see $b \in X \cup Xxy^{-1}.$  By considering $\varphi(\rho^2 \cdot \rho^2 \gamma) \sim b\gamma$ we know $b \in Y \cup Yx^{-1}y.$  Thus $b \in \{1, x^{-1}y^{-1}\}.$  
Assume $b = x^{-1}y^{-1}.$ Since $\rho^2 \gamma = x^{-1}y\gamma \rho^2$ this gives
\[
\rho^2 \gamma = \varphi(x^{-1}y\gamma \cdot \rho^2) \sim x^{-1}y\gamma x^{-1}y^{-1}\rho^2 = x^{-2}y^2 \gamma \rho = x^{-1}y\rho^2 \gamma
\]
which is a contradiction according to Lemma \ref{lemclp4mg}.  Thus $\varphi(\rho^2) = \rho^2.$ 
We may now apply Theorem \ref{thmIdonRot} to obtain  $\varphi = \Id_G$.  Thus the groups of type  {\bf p4mg} have  trivial weak Cayley table group.
\qed\medskip 

Thus we have found generators for $\mathcal W(G)$ for each wallpaper group $G$. Presentations for  $\mathcal W(G)$ can now be obtained; however in the trivial cases one can find such a presentation in \cite {gw}.  We thus only give presentations in the non-trivial cases:

\noindent {\it Groups $G$ of type {\bf p3}}
Here we have 
\begin{align*}
\mathcal{W}(G) &= (\langle \psi_{x}, \psi_{y} \rangle \rtimes \langle \tau, I_r \rangle ) \times \langle \iota \rangle  
 \cong (\mathbb{Z}^2 \rtimes (C_3 \times C_3)) \times C_2.
\end{align*}
Note that $I_x=\psi_x^{-1}\psi_y^{-1}, I_y=\psi_x\psi_y^{-2}$.
Here the actions for the semi-direct product are determined by:
$$\psi_x^{I_r} = I_y^{-1}, \,\,\, \psi_y^{I_r} = \psi_x\psi_y^{-1},   \,\,\,   \tau^{\psi_x} =\tau \psi_x \psi_y,  \,\,\, 
\tau^{\psi_y} = \tau \psi_x^{-1} \psi_y^2.$$
\medskip

\noindent {\it Groups $G$ of type {\bf p4}}
Let $\psi_1: x \mapsto x, \, y \mapsto y^{-1}, \, r \mapsto r^{-1}$. Then we have 
\begin{align*}
\mathcal{W}(G) &= (\langle \tau_{x}, \tau_{y }, \tau_{r^2}, \sigma_x, \sigma_y, \sigma_{r} \rangle \rtimes (\langle I_x, I_y \rangle \rtimes \langle I_r, \psi_1 \rangle )) \times \langle \iota \rangle \\
 &\cong \left (\left [(\mathbb{Z}^2 \rtimes C_2) \times (\mathbb{Z}^2 \rtimes C_4)\right ] \rtimes (\mathbb{Z}^2 \rtimes D_{8})\right ) \times C_2
\end{align*}
Here $\mathbb{Z}^2 \rtimes C_2=\langle \tau_{x}, \tau_{y }, \tau_{r^2}\rangle$ where 
$ \tau_{r^2}^{\psi_x} = \tau_{r^2}\tau_x^{-1}\tau_y,    \tau_{r^2}^{\psi_y} =\tau_{r^2} \tau_x \tau_y^{-1},$ gives this semi-direct product, and
 $\mathbb{Z}^2 \rtimes C_4=\langle\sigma_x, \sigma_y, \sigma_{r} \rangle,$ where
 $  \sigma_x^{\sigma_{r}} = \sigma_y^{-1}  ,\sigma_y^{\sigma_{r}} = \sigma_x$ gives this semi-direct product. 

Next, the subgroup $\mathbb{Z}^2 \rtimes D_{8}$ is 
$\Aut(G)=\langle \psi_x,\psi_y\rangle\rtimes  \langle I_\rho,\psi_1\rangle$. To see that $\mathbb{Z}^2 \rtimes D_{8}$  acts on $\mathbb{Z}^2 \rtimes C_2$ and $\mathbb{Z}^2 \rtimes C_4$ we note that for $\mu  \in \Aut(G)$ we have
$$(\tau_u)^{\mu}=\tau_{\mu^{-1}(u)},\quad (\sigma_u)^{\mu}=\sigma_{\mu^{-1}(u)}.$$ 
Thus for $G$ of type {\bf p4}  we can take $\mathcal N(G)=\langle \tau_{x}, \tau_{y }, \tau_{r^2},  \sigma_x, \sigma_y, \sigma_{r}\rangle$ in Theorem 
\ref {thm1.1}. 
\medskip 

\noindent {\it Groups $G$ of type {\bf p6}} Let $\alpha=xy, \beta=xy^{-2}$ and $\psi: x \mapsto y,  y  \mapsto x,  r  \mapsto  r^{-1}$. 
Here we have 
\begin{align*}
\mathcal{W}(G) &= (\langle \tau_{x^2}, \tau_{y^2}, \tau_{r^3}, \sigma_\alpha, \sigma_\beta, \sigma_{r^2} \rangle \rtimes (\langle I_x, I_y \rangle \rtimes \langle I_r, \psi \rangle )) \times \langle \iota \rangle \\
 &\cong \left (\left [(\mathbb{Z}^2 \rtimes C_2) \times (\mathbb{Z}^2 \rtimes C_3)\right ] \rtimes (\mathbb{Z}^2 \rtimes D_{12})\right ) \times C_2,
\end{align*}
with the details being similar to the last case.\qed\medskip

\noindent {\it Groups $G$ of type {\bf p2mm}} We note that if $G=H \times J$ is a direct product of non-abelin groups, then $G$ always has a non-trivial weak Cayley table isomorphism: just take $\varphi \times \psi$, where $\varphi \in \Aut(H)$ and $\psi \in \mathcal W(J)$ is an anti-automorphism. A group of type  {\bf p2mm} is a direct product. In this case we have
\begin{align*}
\mathcal W(G)=&
\left ( \left (\langle \psi_{1,0}, \psi_{0,1} \rangle \rtimes \langle \tau_s \sigma_{rs}, \rangle \right) \rtimes \langle I_r, I_s, \psi \rangle \right) \times \langle \iota \rangle \\
 &\cong  \left( \left(\mathbb{Z} \rtimes (C_2 \times C_2) \right) \rtimes D_8 \right) \times C_2.\qed
\end{align*}

     \section {The semi direct product cases: $H=\mathbb Z^n \rtimes \mathcal C_2$}

           So let $H=A \rtimes_\theta  \mathcal C_2,$ where $A=\langle a_1,\dots,a_n\rangle\cong \mathbb Z^n$  and $\theta \in \Aut(A).$  Let  $\mathcal C_2=\langle r\rangle$ and $Z=Z(H)$. We write elements of $H$ as $ar^\varepsilon, a \in A,  \varepsilon \in \{0,1\}$. We also assume that $H$ is not abelian. This allows us to assume that $a_i \notin Z(H)$ for all $1\le i\le n$, since if $a_i \in Z(H)$ and $a_j \notin Z(H)$, then we just replace $a_i$ by $a_ia_j$.
           We further note that $a \in A$ is in $Z(H)$ if and only if $a^r=a$.

           Let $A_i=\langle a_i\rangle, 1\le i\le n$.

 \begin {lem}\label{lem2246}  
          (i)  For $a  \in A$ we have $a^{H}=\{a, a^ r\}$. 

(ii)  For $h=ar \in H\setminus A$ we have  $h^H=H'h$.
\end{lem}
\noindent{\it Proof} (i) is clear since $A$ is abelian of index $2$. 

For (ii) we note that 
      $r^a=a^{-1}a^r r,$
     so that for $b \in A$ we have $(br)^{a^m}=b(a^{-1}a^r)^m r$ for all $m \in \mathbb Z$. Let $$\beta_i=a_i^{-1}a_i^r \in H', \quad 1\le i\le n.$$ 
     
     Let $B_i=\langle \beta_i\rangle \le H', 1\le i\le n$, and note that 
     $$( a_i^{-1}a_i^r)^r=(a_i^r)^{-1}a_i=(a_i^{-1}a_i^r)^{-1} \in B_i.$$ 
     
     Also for $a \in A$ we have $(ar)^{A_i}=aB_ir$ and so $(ar)^A=a\langle B_1,\dots,B_n\rangle r$. 
    
     Also 
     \begin{align*}
     (ar)^{Ar}&= (\langle B_1,\dots,B_n\rangle ar)^r=\langle B_1,\dots,B_n\rangle^r a^rr=\langle B_1,\dots,B_n\rangle a^rr\\&
     =   \langle B_1,\dots,B_n\rangle (a^ra^{-1}) \cdot ar.
     \end{align*}     
     
       Here we note that $a^ra^{-1} \in H'$ and the result will follow upon showing that $H' = \langle B_1,\dots,B_n\rangle$.
          Now for $b \in A$ we have
     $$ (a_i,br)=a_i^{-1}rb^{-1}a_ibr=a_i^{-1}ra_ir=a_i^{-1}a_i^r \in B_i,
     $$ and
      $$ (a_ir,a_jr)=ra_i^{-1}ra_j^{-1}a_ira_jr=(a_i^r)^{-1}a_i \cdot a_j^{-1}a_j^r  \in B_iB_j,
     $$ 
          shows that $H' \le \langle B_1,\dots,B_n\rangle$, and it  follows that $H' = \langle B_1,\dots,B_n\rangle$.
     \qed\medskip

\begin {lem}\label{lem2286} Let $\varphi \in \mathcal W(H)$. Then $\varphi(A)=A$ and if $\varphi(a_i)=\alpha_i, 1\le i\le n$, then  $\varphi(a_1^{\lambda_1}\dots a_n^{\lambda_n})=\alpha_1^{\lambda_1}\dots \alpha_n^{\lambda_n}$ for all $ \in \mathbb Z$.
\end{lem}
\noindent{\it Proof} 
          From Lemma \ref {lem2246}  we see that 
          $(ar)^{H}$ is infinite for any $a \in A$. Thus only the elements of $A$ have finite conjugacy classes, and so $\varphi(A)=A$.
Now from $a^H=\{a,a^r\}$ we have  
          $\varphi(a)\sim \varphi(a^r)=\varphi(a)^r$. Thus $\varphi(\{a,a^r\})=\{\varphi(a),\varphi(a)^r\}$.

         \begin {lem}\label{lem2256} For all $n \in \mathbb Z$ and $a \in A$ we have $\varphi(a^n)=\varphi(a)^n$.
         
\end{lem}
         \noindent{\it Proof} We have $a^H=\{a,a^r\}$. 
Let  $\varphi(a)=\alpha,$ so that  $\varphi(a^r)=\alpha^r$.  
          We then have $\varphi(a^\varepsilon)=\alpha^\varepsilon$ for $\varepsilon =0,\pm 1$. So now assume (inductively) that 
          $\varphi(a^i)=\alpha^i$ for $|i|<n$ for some $n\ge 2$ and that $\varphi(a^n)\ne\alpha^n$. Then 
          $$\varphi(a^n)=\varphi(a \cdot a^{n-1})\sim \alpha \alpha^{n-1}\in \{\alpha^n, (\alpha^n)^r\},$$ which shows that $\varphi(a^n)= (\alpha^n)^r$. 
          In particular we have $\varphi(a^n)= (\alpha^n)^r \ne \alpha^n$.
          
          Now 
          $
         \alpha=\varphi(a)=\varphi(a^n\cdot a^{1-n}) \sim (\alpha^r)^n\alpha^{1-n}.
          $
          Thus we must have either 
          
          \noindent (a) $ (\alpha^r)^n\alpha^{1-n}=\alpha$; or 
          
          \noindent (b) 
         $ (\alpha^r)^n\alpha^{1-n}=\alpha^r$.

         If we have (a), then  $ (\alpha^r)^n=\alpha^{n}$, from which we get (since $A$ is a free abelian group)  $\alpha^r=\alpha$, giving a contradiction.
         
         If we have (b), then we get  $(\alpha^r)^{n-1}=\alpha^{n-1}$, which again gives $\alpha^r=\alpha$, since $n\ge 2$. This gives another contradiction and we have completed the proof.\qed\medskip

         For any $a \in A$ we can write $a=a_1^{\lambda_1} a_2^{\lambda_2} \dots a_n^{\lambda_n}, \lambda_i \in \mathbb Z$.  We define a length function on $A$ as follows:
         $$|a|
        =|\lambda_1|+|\lambda_2|+\dots +|\lambda_n|.
        $$ 
          Thus  we now know that $\varphi( a)=\prod_{i=1}^n\alpha_i^{\lambda_i}$ for all $a \in A$ with $|a| \le 1$.

          So now assume (inductively) that there is $m\ge 2$ such that   $\varphi( a_1^{\lambda_1} a_2^{\lambda_2} \dots a_n^{\lambda_n})=\alpha_1^{\lambda_1} \alpha_2^{\lambda_2} \dots \alpha_n^{\lambda_n}$ for all $a_1^{\lambda_1} a_2^{\lambda_2} \dots a_n^{\lambda_n} \in A$ with $|a_1^{\lambda_1} a_2^{\lambda_2} \dots a_n^{\lambda_n}| <m$, and that there is $a_1^{\lambda_1} a_2^{\lambda_2} \dots a_n^{\lambda_n}$ such that
      $|a_1^{\lambda_1} a_2^{\lambda_2} \dots a_n^{\lambda_n}| =m$ and       $\varphi(a_1^{\lambda_1} a_2^{\lambda_2} \dots a_n^{\lambda_n})\ne \alpha_1^{\lambda_1} \alpha_2^{\lambda_2} \dots \alpha_n^{\lambda_n}$. 
      Without loss we can assume that $\lambda_1\ge 1$.

      Further, since  $\varphi(a_1^{\lambda_1} a_2^{\lambda_2} \dots a_n^{\lambda_n})\ne \alpha_1^{\lambda_1} \alpha_2^{\lambda_2} \dots \alpha_n^{\lambda_n}$, we must have
  $\varphi(a_1^{\lambda_1} a_2^{\lambda_2} \dots a_n^{\lambda_n})= (\alpha_1^{\lambda_1} \alpha_2^{\lambda_2} \dots \alpha_n^{\lambda_n})^r$     
since
          \begin{align*}
            \varphi( a_1^{\lambda_1} a_2^{\lambda_2} \dots a_n^{\lambda_n}  )&=\varphi(a_1 \cdot a_1^{\lambda_1-1} a_2^{\lambda_2} \dots a_n^{\lambda_n})
            \sim \varphi(a_1)\varphi(a_1^{\lambda_1-1} a_2^{\lambda_2} \dots a_n^{\lambda_n})\\&=\alpha_1 \cdot \alpha_1^{\lambda_1-1} \alpha_2^{\lambda_2} \dots \alpha_n^{\lambda_n}=\alpha_1^{\lambda_1} \alpha_2^{\lambda_2} \dots \alpha_n^{\lambda_n}.
            \end{align*}

      By induction we   have $\varphi(   a_2^{\lambda_2} \dots a_n^{\lambda_n}  )= \alpha_2^{\lambda_2} \dots \alpha_n^{\lambda_n}$; however,  using Lemma \ref {lem2256}, we also have
      $$\varphi( a_2^{\lambda_2} \dots a_n^{\lambda_n})
      =\varphi( a_1^{\lambda_1}  a_2^{\lambda_2} \dots a_n^{\lambda_n}\cdot a_1^{-\lambda_1})\sim 
      (\alpha_1^{\lambda_1} \alpha_2^{\lambda_2} \dots \alpha_n^{\lambda_n})^r \cdot \alpha_1^{-\lambda_1}
      $$
      Thus we have
      $\alpha_2^{\lambda_2} \dots \alpha_n^{\lambda_n} \sim (\alpha_1^{\lambda_1} \alpha_2^{\lambda_2} \dots \alpha_n^{\lambda_n})^r \cdot \alpha_1^{-\lambda_1}.
      $

         This gives two cases:
      
 \noindent (a)     $\alpha_2^{\lambda_2} \dots \alpha_n^{\lambda_n} = (\alpha_1^{\lambda_1} \alpha_2^{\lambda_2} \dots \alpha_n^{\lambda_n})^r \cdot \alpha_1^{-\lambda_1};$ or 
 
  \noindent (b)     $\alpha_2^{\lambda_2} \dots \alpha_n^{\lambda_n} = (\alpha_1^{\lambda_1} \alpha_2^{\lambda_2} \dots \alpha_n^{\lambda_n}) \cdot (\alpha_1^{-\lambda_1})^r.$ 
  
  \medskip

       Here possibility  (a) gives $(\alpha_1^{\lambda_1} \alpha_2^{\lambda_2} \dots \alpha_n^{\lambda_n})^r =\alpha_1^{\lambda_1} \alpha_2^{\lambda_2} \dots \alpha_n^{\lambda_n}$, a contradiction since $\varphi(a_1^{\lambda_1} a_2^{\lambda_2} \dots a_n^{\lambda_n})=(\alpha_1^{\lambda_1} \alpha_2^{\lambda_2} \dots \alpha_n^{\lambda_n})^r \ne \alpha_1^{\lambda_1} \alpha_2^{\lambda_2} \dots \alpha_n^{\lambda_n}$. 
        
        But (b) gives  $ (\alpha_1^{\lambda_1})^r=\alpha_1^{\lambda_1}$; this shows that $\varphi(a_1)=\alpha_1 \in Z(H)$, which in turn shows that $a_1 \in Z(H)$, 
        which is  a contradiction to the choice of the $a_i$. This concludes the proof of Lemma  \ref {lem2286}. \qed\medskip

        Thus we see  that $\varphi|_A:A \to A$ is an automorphism of $A$.

          \begin {lem}\label{lemaut25889} 
 (1) For any $b \in A$  such that $b r$ has order $2$ the homomorphism determined by $a \mapsto a, (a \in A), r \mapsto b r$ is   an automorphism of $H$.
 
 (2)  Let $\psi \in \Aut(A)$ where $\psi$ commutes with $\theta$. Then the homomorphism determined by $a \mapsto \psi(a), (a \in A), r \mapsto  r$ is   an automorphism of $H$.
  \end{lem}
\noindent{\it Proof} (1) is clear since $r$ and $b r$ have order $2$ and $a^{b r}=a^r$.

(2) The relations in $H$ have the form $a^r=a^\theta$ and $r^2=1$. We have
$\psi(a^r)=\psi(\theta(a))=\theta(\psi(a))=(\psi(a))^r,
$
and we are done.
\qed\medskip

    Since $\varphi(Ar)=Ar$ we can assume that $\varphi(r)=b r$ for some $b \in A$.  Further, we know that $b r$ has order $2$, since $\varphi$ is a weak Cayley table map.

           By Lemma \ref {lemaut25889}  we can assume (by composing $\varphi$ with the inverse of such an automorphism)  that $b=1$ i.e. that $\varphi(r)=r$.
           
                     \begin {lem}\label{lem32581} Suppose that 
 $a \in A$. Then $\varphi(a^r)=\varphi(a)^r$.
 \end{lem}
\noindent{\it Proof} We know that  $\varphi(a^r) \in \{\varphi(a), \varphi(a)^r\}$; however if $\varphi(a^r) =\varphi(a)$, then $a^r=a$ (so that $a=a^r$ and $\varphi(a)=\varphi(a^r)$ are central), and otherwise we have $\varphi(a^r) =\varphi(a)^r$. 
           \qed\medskip
           
           This lemma shows that the automorphism of $A$ determined by the action of $\varphi$ commutes with the action of $\theta$. 
           Since $\varphi|_A$ is a automorphism of $A$ we can now use Lemma   \ref {lemaut25889}  (2) and Lemma \ref {lem32581}  so as to be able to assume that $\varphi|_A=id|_A$. 
           
           Thus we can now assume: $\varphi|_A=Id|_A$ and $\varphi(r)=r$.

          By Lemma \ref {lemra=rb} we know that if 
 $a \in A$ and  $\varphi(ar)=br, b \in A$, then $b\in \{a,a^r\}$.            
           Thus  if $\varphi(ar) \ne ar$, then   $\varphi(ar)=a^{r}r=ra$.  
              
           Now assume that there are $a,b \in A \setminus \{1\}$ such that $\varphi(ar)=ar\ne a^r r, \varphi(br)=b^{r}r \ne br$. Then
           $$ba^{r}=\varphi(ba^{r})=\varphi(br\cdot ar)\sim \varphi(br)\varphi(ar)=b^{r}rar=b^{r}a^{r}.$$
           
           It follows that either $ba^{r}=b^{r}a^{r}$ or $ba^{r}=ab$, each possibility giving a contradiction.

           Thus  there is $\varepsilon \in \{0,1\}$ such that $\varphi(ar)=a^{r^\varepsilon} r$ for all $a \in A$. If $\varepsilon=0$, then $\varphi$ is the identity map on $H$. 
          
           So now assume that  $\varepsilon=1$, so that   $\varphi$ satisfies
           $$\varphi(a)=a,\quad \varphi(ar)=a^{r}r,\quad \text { for all } a \in A.$$
    Then for all $a,b \in A$ we have:
 $ rab=(ab)^rr=\varphi(a\cdot br) \sim a b^rr=arb,$ so that either (i) $rab=arb$; or (ii) $rab=(arb)^r$.  But (i) gives $(a,r)=1$, and (ii) gives $(b,r)=1$. Since (i) and (ii) are true for all $a,b \in A$ we have a contradiction,
    and we have proved
    Theorem \ref {thmh}.\qed\medskip


                          \begin {thm}\label{lemgex} 
               Let  $p$ be an odd  prime, $A$  an abelian group and let $G= A \rtimes_\theta  \mathcal C_p$ be a semi direct product.
We   assume that $\theta$ is not trivial. Then $G$ has a non-trivial weak Cayley table map.
  \end{thm}
\noindent{\it Proof} This will also give a proof of Theorem \ref {thmh2}.   
                   Let $\mathcal C_p=\langle r\rangle.$
                  
                           \begin {lem}\label{lemgex} 
         (i)        Let $1\le k<p$.           Then every $g \in G'$ has the form $(a,r^k)$ for some $a \in A$.

(ii) For any $1\le k <p$ and $g \in G'$ there is $a \in A$ such that $(r^k)^a=gr^k.$
  \end{lem}
\noindent{\it Proof} The proof of  (i) is standard, using the Witt-Hall identities \cite [ p. 290] {mks}. 


               
               
             (ii)  Now let $g \in G'$ and fix $1 \le k<p$. Then by  (i)  there is $a \in A$ such that 
               $$g=(a,r^{-k})=a^{-1}r^kar^{-k}, \text { so that } gr^k=(r^k)^a.
               $$
               This concludes the proof of Lemma \ref {lemgex}.\qed\medskip

Now define a bijection $\varphi: G \to G$ by
$$\varphi(a)=a^r \text { for } a \in A, \quad \varphi(g)=g \text { for } g \in G \setminus A.$$

Since $\theta$ is not trivial there is some $a \in A$ such that $a^r \ne a$. 
Thus $ \varphi(a)\varphi(r)=a^r \cdot r \ne ar =\varphi(a r),$ so  that $\varphi$ is not a homomorphism. Similarly we have, using the fact that $p>2$, 
 $$ \varphi(ar)\varphi(r)=arr \ne rar =\varphi(r \cdot ar),$$ so  that $\varphi$ is not an antihomomorphism. 

We now show that $\varphi$ is a weak Cayley table isomorphism. We need to show that $\varphi(g_1g_2) \sim \varphi(g_1)\varphi(g_2)$ for all $g_1,g_2 \in G$, and we do this by considering the following situations:

\noindent {\bf CASE 1:} $g_1,g_2 \in A$. Here we have $\varphi(g_1) \varphi(g_2)=g_1^rg_2^r=(g_1g_2)^r=\varphi(g_1g_2).$
\medskip

   \noindent {\bf CASE 2:} $g_1\in A,g_2 \notin A$. 
   Here we write $g_1=a, g_2=br^k, a,b \in k, 1\le k<p$. 
Then $$\varphi(g_1g_2)=\varphi(abr^k)=abr^k, \quad \varphi(g_1)\varphi(g_2)=\varphi(a)\varphi(br^k)=(a^rb)\cdot r^k.
$$
Now by Lemma \ref{lemgex} there is $c \in A$ such that 
$$(a^rb)(ab)^{-1}=a^ra^{-1}=(c,r^{-k})=(r^k)^c\cdot r^{-k}.
$$ 
Then we have
\begin{align*}
\varphi(g_1g_2)^c&=(abr^k)^c=(ab)(r^k)^c=(ab)\cdot (a^rb)(ab)^{-1}r^k\\&
=a^r(br^k)=  \varphi(g_1)\varphi(g_2).
\end{align*}

      This does this case and
     a similar argument does the case  where   $g_1\notin A, g_2 \in A$.
   \medskip

    \noindent {\bf CASE 3:} $g_1\notin A, g_2 \notin A, g_1g_2 \notin A$. Here we have  $\varphi(g_1)=g_1, \varphi(g_2)=g_2, \varphi(g_1g_2)=g_1g_2$ and so     $\varphi(g_1g_2)=g_1g_2=\varphi(g_1)\varphi(g_2),$
   which does this case.
   \medskip

    \noindent {\bf CASE 4:} $g_1\notin A, g_2 \notin A, g_1g_2 \in A$.   Here we have
     $\varphi(g_1g_2)=(g_1g_2)^r$ and $\varphi(g_1) \varphi(g_2)=g_1g_2$, from which it follows that 
     $\varphi(g_1g_2)=(\varphi(g_1) \varphi(g_2))^r$.   \medskip
             
      This concludes the proof that $\varphi$ is a weak Cayley table isomorphism.
      \qed


\begin{thebibliography}{JMS}
 








\bibitem[Cur]{cur} Curtis, Charles W., \emph{Pioneers of
representation theory: Frobenius, Burnside, Schur, and Brauer}. History of Mathematics {\bf 15}, American
Mathematical Society, Providence, RI; London Mathematical Society, London, (1999), 287 pages.

\bibitem[CR]{cr} Curtis, Charles W.; Reiner, Irving, \emph{ Representation theory of finite groups and associative algebras.} Reprint of the 1962 original. AMS Chelsea Publishing, Providence, RI, (2006).




\bibitem[Fe]{fe} Feit W., Characters of finite groups, W. A. Benjamin, (1967).



 
 



\bibitem[GW]{gw}    Gon\c{c}alves, Daciberg; Wong, Peter \emph{Automorphisms of the two dimensional crystallographic groups}. Comm. Algebra 42 (2014), no. 2, 909--931.
 
 














 \bibitem[Hu]{hu} Humphries, Stephen P., \emph{Weak Cayley table groups}, J. Algebra {\bf 216} (1999), 135-158.



 \bibitem[HN]{hu2} Humphries, Stephen P., Nguyen Long, \emph{Weak Cayley table groups: Alternating and Coxeter groups}, preprint (2013), 16 pages.


 \bibitem[HN2]{hu3} \bysame, \emph{Weak Cayley Table groups III: $PSL(2,q)$}, preprint (2014). 30 pages.














\bibitem[Iv]{iv}   Iversen, Birger \emph{Lectures on crystallographic groups}. Lecture Notes Series, 60. Aarhus Universitet, Matematisk Institut, Aarhus, (1990) vi+144 pp. 
 


\bibitem[Ja]{ja}  Janssen, T.
\emph{Crystallographic groups}. North-Holland Publishing Co., Amsterdam-London; American Elsevier Publishing Co., Inc., New York, (1973) xiii+281 pp. 


 \bibitem[JMS]{jms} Johnson, Kenneth W.; Mattarei, Sandro; Sehgal, Surinder K., \emph{ Weak Cayley tables}, J. London Math. Soc., {\bf 61} (2000), 395--411.
















\bibitem[MA]{ma} Bosma W., Cannon J., MAGMA (University of Sydney), (1994). 

\bibitem[MKS]{mks} Magnus, W, Karrass, Solitar, \emph{Combinatorial group theory}, Dover, (1976)





























  
 \end{thebibliography}
\end{document}